\documentclass[twoside,11pt]{article}
\usepackage{journal2e_v1}
\usepackage{epsf}
\usepackage{epsfig}
\usepackage{amsmath}
\usepackage{subfigure}
\usepackage{amssymb}
\usepackage{multirow}
\usepackage{multicol}
\usepackage{diagbox}
\usepackage{makecell}

\usepackage{graphicx} 

\usepackage{booktabs}
\usepackage{algorithm}
\usepackage{algorithmic}
\usepackage{color}
\usepackage[pagebackref=true,breaklinks=true,letterpaper=true,colorlinks,bookmarks=yes]{hyperref}

\newtheorem{assumption}[theorem]{Assumption}

\def\A{{\bf A}}

\def\g{{\bf g}}

\def\u{{\bf u}}

\def\v{{\bf v}}

\def\w{{\bf w}}
\def\X{{\bf X}}
\def\x{{\bf x}}

\def\y{{\bf y}}

\def\0{{\bf 0}}
\def\1{{\bf 1}}

\def\denote{\stackrel{\Delta}{=}}

\def\AM{{\mathcal A}}

\def\RB{{\mathbb R}}

\def\NB{{\mathbb N}}

\def\argmin{\mathop{\rm argmin}}

\def\dom{\mathrm{dom}}

\def\sgn{\mathrm{sgn}}

\def\st{\mathsf{s.t.}}

\begin{document}
\title{On the Global Convergence of Majorization Minimization Algorithms for Nonconvex Optimization Problems}

\author{\name Yangyang Kang  \\
\addr Department of Computer Science and Engineering \\
Shanghai Jiao Tong University \\ 
\texttt{kcany27@gmail.com}
\AND
\name Zhihua Zhang \\
\addr Department of Computer Science and Engineering \\
Shanghai Jiao Tong University \\
\texttt{zhang-zh@cs.sjtu.edu.cn}
\AND
\name Wu-Jun Li \\
\addr National Key Laboratory for Novel Software Technology \\
Department of Computer Science and Technology \\
Nanjing University \\
\texttt{liwujun@nju.edu.cn}
}
\date{}

\maketitle

\begin{abstract}%
In this paper, we study the global convergence of majorization minimization (MM) algorithms
for solving nonconvex regularized optimization problems. MM algorithms  have received great
attention in machine learning. However, when applied to nonconvex optimization problems,
the  convergence of MM algorithms is a challenging issue.
We introduce theory of the Kurdyka-\L{}ojasiewicz inequality to address this issue. In particular,
we show that many nonconvex problems enjoy the  Kurdyka-\L{}ojasiewicz property and establish the global convergence result of
the corresponding MM procedure. We also extend  our result to a well known method that called CCCP (concave-convex procedure).
\end{abstract}
\section{Introduction}

Majorization minimization (MM) algorithms  have
wide applications in machine learning and statistical inference~\citep{lange2000optimization,lange2004optimization}. The MM algorithm can be regarded as a generalization of expectation-maximization
(EM) algorithms, and it aims to turn an otherwise hard or complicated optimization problem into a tractable one by alternatively iterating an \emph{Majorization} step and an  \emph{Minimization} step.

More specifically, the majorization step constructs a tractable surrogate function to substitute the original objective function
and the minimization step minimizes this surrogate function to obtain a new estimate of  parameters in question.
In the conventional MM algorithm,
convexity plays a key role in the construction of surrogate functions. Moveover, convexity arguments make  the conventional MM algorithm have
the same convergence properties as EM algorithms~\citep{lange2004optimization}.

Alternatively, we are interested in  use of MM algorithms in solving nonconvex (nonsmooth) optimization problems.
For example,  nonconvex penalization has been demonstrated to have attractive properties  in sparse estimation.
In particular,
there exist many nonconvex penalties, including the $\ell_q$ ($q\in (0,1)$) penalty,  the
smoothly clipped absolute deviation (SCAD) \citep{fan2001variable},
the minimax concave plus penalty (MCP) \citep{zhang2010nearly},  the capped-$\ell_1$ function~\citep{zhang2010analysis,zhang2012general,gong2013general},  
the LOG penalty~\citep{mazumder2011sparsenet,armagan2013generalized},  etc.
However, they might yield computational challenges due to  nondifferentiability and nonconvexity that they have.
An MM algorithm would be a desirable choice~\citep{lange2004optimization}.

In this paper we would like to address the global convergence property of MM algorithms for nonconvex optimization problems.
Our  motivation comes from  the novel Kurdyka-\L{}ojasiewicz inequality. In the pioneer work \citep{lojasiewicz1963propriete,lojasiewicz1993geometrie}, the author provided the ``\L{}ojasiewicz inequality" to derive  finite trajectories.
Later on,  \cite{kurdyka1998gradients} extended the \L{}ojasiewicz inequality to definable functions and applications.
\cite{bolte2007lojasiewicz} then extended to nonsmooth subanalytic functions.
Recently,  the Kurdyka-\L{}ojasiewicz property
has been used to establish convergence analysis  of proximal alternating minimization or coordinate descent algorithms~\citep{attouch2010proximal,xu2013block,bolte2013proximal}.

We revisit  a generic MM procedure of solving  nonconvex optimization problems.
We observe that many nonconvex penalty functions satisfy the Kurdyka-\L{}ojasiewicz inequality and
such a property is shared by a number of machine learning problems arising in a wide variety of applications. Specifically, we demonstrate several examples, which admit the
Kurdyka-\L{}ojasiewicz  property.
Thus, we conduct the convergence analysis of the MM procedure based on theory of the Kurdyka-\L{}ojasiewicz inequality.
More specifically, our work offers the following major contributions.
\begin{itemize}

\item
We discuss a family of nonconvex optimization problems in which the objective function consists of a smooth function and a non-smooth function.
We  give the constructive criteria of surrogates that approximate the original functions well.  Additionally,
we also illustrate that many existing methods for solving the nonconvex optimization problem can be regarded as an MM procedure.
\item  We establish the global convergence results of a generic MM framework for the nonconvex problem  which are obtained  by exploiting the geometrical property of the objective function around its critical point. To the best of our knowledge, our work is the first study to address  the convergence property of MM algorithms for nonconvex optimization using the Kurdyka-\L{}ojasiewicz  inequality.
\item We also show that our global convergence results can be successfully extended to many popular and powerful methods such as
iteratively re-weighted $\ell_1$ minimization method~\cite{candes2008enhancing,chartrand2008iteratively},  local linear approximation (LLA) \cite{zou2008one,zhang2010analysis},  concave-convex procedure (CCCP)~\cite{yuille2003concave,lanckriet2009convergence}, etc.
\end{itemize}

\subsection{Related Work and Organization}\label{sec:related}

We  discuss some related work about the convergence analysis of nonconvex
optimization. \cite{vaida2005parameter} established the global convergence of EM algorithms and  extended it to the global convergence of MM algorithms under some conditions. However, they considered the differentiable objective function, whereas the objective function in our paper can be nonsmooth (also nonconvex). This implies that the problem we are considering is more challenging. Additionally, \cite{vaida2005parameter} assumed that all the stationary points of objective function are isolated. In our paper, we don't require this assumption. The isolation assumption does not always hold, or holds but is difficult to verify, for many objective functions in practice. This motivates us to employ the Kurdyka-Lojasiewicz inequality to establish the convergence. Moreover, it is usually easily verified that the objective function admits the Kurdyka-Lojasiewicz inequality. \cite{gong2013general}  proposed an efficient iterative shrinkage and thresholding algorithm  to solve nonconvex regularized problems.
The key assumption is that the computation of proximal operator of the regularizer has a closed form. We note that this method
falls into our MM framework. However, the authors only showed that the subsequence converges to a critical point.
\cite{mairal2013optimization} studied instead asymptotic stationary
point conditions with first-order surrogate functions, but  he did not propose the convergent sequence which converges to the solution point.

\cite{attouch2010proximal,xu2013block,bolte2013proximal} employed the Kurdyka-\L{}ojasiewicz inequality to analyze the convergence of
nonconvex optimization problems. They are mainly concerned with the convergence analysis of the block coordinate approaches. In this paper, we pay attention to the global convergence analysis of the MM framework for solving nonconvex regularization problems. Specifically, we construct surrogates both on the smooth and nonsmooth terms. To achieve the global convergence, we exploit the geometry property of the objective function around its critical point.

The remainder of the paper is organized as follows. Section \ref{sec:pre} provides preliminaries about the nonsmooth and nonconvex analysis and introduces the Kurdyka-\L{}ojasiewicz property. We also give some examples which enjoy the Kurdyka-\L{}ojasiewicz inequality. In Section \ref{sec:problem}, we formulate the problem we are interested in and make some common assumptions. A generic majorization minimization algorithm is revisited in Section \ref{sec:MM}. Section \ref{sec:Convergence} is the key part of our paper which gives the  global convergence results.
In Section \ref{sec:cccp} we extend our work to CCCP. In Section \ref{sec:experiment} we conduct numerical examples to verify our theoretical results. Finally, we conclude our work in Section \ref{sec:conclusion}.

\section{Preliminaries}
\label{sec:pre}

In this section we introduce the notion of Fr\'{e}chet's subdifferential and a limiting-subdifferential.
Then we present the novel  Kurdyka-\L{}ojasiewicz inequality.
 First of all,  for any $\u=(u_1, \ldots, u_p)^T\in\RB^p$ and $\v=(v_1, \ldots, v_p)^T \in\RB^p$,  we denote $\langle\u,\v\rangle=\sum^p_{i=1}u_iv_i$ and $\|\u\|=\sqrt{\langle\u,\u\rangle}$
 here and later.

\begin{definition}[Subdifferentials] \citep{rockafellar1998variational}
\label{def:subdiff}
Consider a proper and lower semi-continuous function $f: \RB^{p} \to(-\infty, +\infty]$ and a point ${\x} \in \dom(f)$.
\begin{enumerate}
\item[\emph{(i)}] The Fr\'{e}chet subdifferential of $f$ at ${\x}$, denoted $\hat{\partial} f({\x})$, is the set of all vectors $\u \in \RB^p$
which satisfy
\[
\liminf_{\begin{footnotesize}\begin{array}{c} \y \neq \x \\ \y \to \x \end{array} \end{footnotesize}}
\frac{f(\y)- f(\x)- \u^T(\y-\x) }{\|\y-\x \|} \geq 0.
\]
\item[\emph{(ii)}] The limiting-subdifferential of $f$ at ${\x}$, denoted ${\partial} f({\x})$, is defined as
\[
{\partial} f({\x}) \equiv \bigg\{\u \in \RB^p: \exists \x_k \to \x, f(\x_k) \to f(\x) \; \mbox{ and } \; \u_k \in \hat{\partial} f({\x_k})
\to \u \; \mbox{ as } \; k\to \infty  \bigg\}.
\]
\end{enumerate}
\end{definition}



\begin{remark}\label{rem:1} Here
$\dom f \denote \big\{\x: f(\x) < +\infty\big\}$. If $\x \notin \dom f$, one sets $\hat \partial f(\x)= \varnothing$.  It is worth pointing out  that $\hat{\partial} f(\x)$ for each $\x$  is closed and convex while ${\partial} f(\x)$ is closed. If $f$ is differentiable at $\x_0$, then $\hat \partial f(\x_0)=\{\nabla f(\x_0)\}$ and $\nabla f(\x_0)\in\partial f(\x_0)$. More details are referred to \cite{rockafellar1998variational}.
As we see, both the Fr\'{e}chet subdifferential and limiting-subdifferential are applicable for nonconvex functions.


\begin{corollary}[Rockafellar and Wets, 1998]\label{def:sum-subdiff}
Suppose $F=f+r:\RB^p\to \RB$. Moreover, $f$ is smooth in the neighborhood of $\x_0$ and $r$ is finite at $\x_0$. Then, we have
\[
\hat\partial F(\x_0)=\nabla f(\x_0)+\hat\partial r(\x_0) \; \mbox{ and } \; \partial F(\x_0)=\nabla f(\x_0)+ {\partial} r({\x_0}).
\]
\end{corollary}

\end{remark}
\begin{definition}\label{def:critical}
It is said that  $\x^* \in \RB^p$ is a critical point of a lower semi-continuous function $F : \RB^p \to \RB \cup \{+\infty\}$,
if  the following condition holds
\[
 \0 \in \partial F(\x^*).
\]
\end{definition}
\begin{remark}\label{rem:2}
If $\x^*$ is a minimizer (not necessarily global) of function F, we can conclude that $\0\in\partial F(\x^*)$.
The set of critical points of $F$ is denoted by $\mathrm{crit} F$.
\end{remark}

\subsection{Kurdyka-\L{}ojasiewicz properties}

With the notion of subdifferentials, we now briefly recall the Kurdyka-\L{}ojasiewicz inequality, which plays a central role in our globally convergence analysis.
\begin{definition}
Let the function $F$  $: \RB^{p}\to(-\infty, +\infty]$ be proper and lower semi-continuous. Then $F$ is said to have the Kurdyka-\L{}ojasiewicz property at $\bar{\u} \in \dom \, \partial F$ if there exist $\eta \in (0, +\infty]$, a neighborhood $\mathcal{U}$ of $\bar{\u}$, and a continuous concave function $\phi:[0, \eta)\to \RB_{+}$ with the following properties:
\begin{enumerate}
\item[\emph{(a)}] $\phi(0)=0$,
\item[\emph{(b)}] $\phi$ is $C^{1}$ on $(0,\eta)$,
\item[\emph{(c)}] for all $t\in(0,\eta)$, $\phi^{\prime}(t)>0$,
\end{enumerate}
such that for all $\u$ in $\mathcal U$ $\cap [F(\bar \u )< F(\u) < F(\bar \u)+\eta]$,  the following {\bf Kurdyka-\L{}ojasiewicz inequality} holds true:
\[\phi^{\prime}(F(\u)- F(\bar \u)) \mathrm{dist}(\0, \partial F(\u))\geq 1.
\]
Here $\mathrm{dist}(\u, \AM)=\inf_\v \left\{\|\u-\v\|, \v\in \AM\right\}$.
\end{definition}

It is well established that real analytic  and  sub-analytic functions satisfy the Kurdyka-\L{}ojasiewicz property~\citep{bolte2007lojasiewicz}.
Moreover, the sum of a real analytic
function and a subanalytic function is  subanalytic \citep{bochnak1998real}. Thus, the sum admits  the Kurdyka-\L{}ojasiewicz property.
Many functions involved in machine learning satisfy the Kurdyka-\L{}ojasiewicz property.
For example, both the logistic loss and the least squares loss are real analytic.

We also find that many nonconvex penalty functions, such as
MCP, LOG, SCAD, and Capped $\ell_1$,  enjoy the Kurdyka-\L{}ojasiewicz property. Here we give two examples. First,
the MCP function is defined as
\[
\zeta(t; \lambda, \gamma)= \left\{ \begin{array}{ll} \lambda(|t|-\frac{t^2}{2\lambda\gamma}) & \mbox{ if } {|t|<\lambda\gamma}, \\
 \frac{\lambda^2\gamma}{2} & \mbox{ if } {|t|\geq\lambda\gamma}, \end{array} \right.
\]
where $\lambda, \gamma>0$ are constants. 
The graph of $\zeta$ is the closure of the following set
\[
\begin{split}
& \bigg\{(t,s):s=\frac{\lambda^2\gamma}{2},t<-\lambda\gamma\bigg\} \cup\bigg\{(t,s):s=\frac{\lambda^2\gamma}{2},t>\lambda\gamma\bigg\} \\
& \cup \bigg\{(t,s):s=-\lambda t-\frac{t^2}{2\gamma},-\lambda \gamma< t<0\bigg\}\cup\bigg\{(t,s):s=\lambda t-\frac{t^2}{2\gamma},
0<t<\lambda\gamma\bigg\}.
\end{split}
\]
This implies that  MCP is a semi-algebraic function \citep{bochnak1998real}, which is  sub-analytic~\citep{bolte2007lojasiewicz}. Thus, MCP satisfies the Kurdyka-\L{}ojasiewicz property. Similarly, we can obtain that the SCAD and capped $\ell_1$ penalties have  the Kurdyka-\L{}ojasiewicz property.

Second, the LOG penalty
is defined as
\[
\zeta(t;\lambda,\alpha)=\lambda\log(1+\alpha|t|), \; \mbox{for } \alpha>0.
\]
The graph of $\zeta$ is the closure of the following set
\[
\begin{split}
&\bigg\{(t,s):s=\lambda\log(1+\alpha t),t>0\bigg\}\cup\bigg\{(t,s):s=\lambda\log(1-\alpha t), t<0\bigg\}.
\end{split}
\]
Note that  the graph is sub-analytic \citep{bolte2007lojasiewicz}, so
the LOG penalty is sub-analytic, which enjoys the Kurdyka-\L{}ojasiewicz property.



\section{Problem and Assumptions}
\label{sec:problem}

In this paper we are mainly concerned with the following optimization problem
\begin{equation}\label{eq:problem}
\min_{\w\in\RB^p} \Big\{F(\w)=f(\w)+r(\w)\Big\}.
\end{equation}
Many machine learning problems can be cast into this formulation.
Typically, $f(\w)$ is defined as a loss function  and $r(\w)$
is defined as a regularization (or penalization) term. 
Specifically, given a training dataset $\mathcal{D}=\big\{(\x_1, y_1),(\x_2, y_2), \cdots, (\x_n, y_n)\big\}$,
one defines $f(\w)=\frac{1}{n}\sum_{i=1}^n f_{i}(\w;\x_i, y_i)$.
A very common setting for the penalty function $r(\w)$ is $\sum_{i=1}^{p}r_i(w_i; \lambda)$,
where $\lambda$ is the tuning parameter controlling the trade-off between the loss function and the regularization.

Recently, many nonconvex penalty functions, such as LOG  \citep{mazumder2011sparsenet,armagan2013generalized}, SCAD~\citep{fan2001variable},  MCP \citep{zhang2010nearly}, and the capped-$\ell_1$ function~\citep{zhang2010analysis},   have been proposed to model sparsity.
These penalty functions have been demonstrated to have attractive properties theoretically and practically.

Meanwhile,  iteratively reweighted methods haven widely used to solve the optimization problem in~(\ref{eq:problem}).
Usually, the iteratively reweighted method enjoys a majorization minimization (MM) procedure.
In this paper we attempt to conduct convergence analysis of the MM procedure. For our purpose, we make some assumptions
about the objective function.

\begin{assumption} \label{ASS1}
Suppose $f : \RB^p \to \RB_+$ is a smooth function of the type $C^{1,1}$. Moreover,  the gradient of $f$ is $L_{f}$-Lipschitz continuous; that is,
\begin{equation}\label{A1}
\|\nabla f(\u)-\nabla f(\v)\|\leq L_{f} \|\u-\v\|
\end{equation}
for any $\u, \v \in \RB^p$, where $L_{f}>0$ is called a Lipschitz constant of $\nabla f$.
\end{assumption}

\begin{corollary}\label{cor:1} Let $h(\w)=\sum_{i=1}^n\alpha_i f_i(\w)$.
Suppose $f_i(\w)$ is differentiable for any $i\in [n] \triangleq \{1, 2, \cdots, n\}$. If each $\nabla f_i(\w)$ is $L_i$-Lipschitz continuous $(L_i>0)$,  then $h(\w)$ is differentiable and $\nabla h(\w)$ is $\sum_{i=1}^n |\alpha_i| L_i$-Lipschitz continuous.
\end{corollary}

\begin{lemma}\label{lemm:1}
If $f : \RB^p \to \RB$ is differentiable and $\nabla f$ is $L_{f}$-Lipschitz continuous. Then
\begin{equation}\label{eq:2}
f(\u)\leq f(\v)+ \langle \nabla f(\v), \u-\v\rangle+\frac{L_{f}}{2}\|\u-\v\|^2
\end{equation}
for any $\u, \v \in \RB^p$.
\end{lemma}
This is a classical result whose proof can be seen from \cite{nesterov2004introductory}.

\begin{assumption}\label{ASS3}
$F : \RB^p \to \RB$ is lower semi-continuous and coercive \footnote {A function $g(\u)$ on $\RB^p$ is said to be coercive if $\lim_{\|\u\|\to \infty} g(\u)=\infty$}, and it satisfies  $\inf_{\w\in\RB^p} F(\w)>-\infty$.
\end{assumption}


We give several examples to show that the assumptions  hold in
many machine learning problems.
For the linear regression, $f(\w)=\frac{1}{2n}\|\X\w-\y\|^2$, where $\w\in\RB^p, \X=[\x_1, \cdots, \x_n]^T \in \RB^{n\times p}$ is the input matrix and $\y=[y_1,\cdots,y_n]^T\in\RB^n$ is the output vector. In this example, the Lipschitz constant of $\nabla f(\w)$ is lower-bounded by the maximum eigenvalue of $\frac{1}{n}\X^T\X$. In binary classification problems in which $y_i \in \{-1, 1\}$,  we consider the logistic regression loss function. Specifically, $f(\w)=\frac{1}{n}\sum_{i=1}^n\log(1+\exp(-y_i\x_i^T\w))$.
The Lipschitz constant of $\nabla f(\w)$ is lower-bounded by $\frac{1}{4n}\sum_{i=1}^n\x_i^T\x_i$.


\begin{table}[t]
\caption{Examples of nonconvex penalties for one dimension}
\label{sample-table}
\vskip 0.15in
\begin{center}
\begin{small}
\begin{sc}
\begin{tabular*}{0.95\textwidth}{l |c}
\hline
\diagbox{Function}& $\zeta(t)$\\
\hline
LOG    & $\frac{\lambda}{\log(\theta+1)}\log(1+\theta|t|)$, $(\theta>0)$ \\
\hline
SCAD & $\left\{ \begin{array}{ll}  \lambda|t| & \mbox{if } {|t|\leq \lambda}, \\
 -\frac{|t|^2-2\theta\lambda|t|+\lambda^2}{2(\theta - 1)} & \mbox{if } {\lambda<|t|\leq\theta\lambda}, \\
  \frac{(\theta+1)\lambda^2}{2} & \mbox{if } \theta \lambda<|t|,  \end{array} \right. (\theta>2)$ \\
\hline
MCP     & $ \left\{ \begin{array}{ll} \lambda(|t|-\frac{t^2}{2\lambda\gamma}) & \mbox{if }  |t|<\lambda\gamma,  \\ \frac{\lambda^2\gamma}{2} & \mbox{if } |t| \geq \lambda\gamma. \end{array} \right. $ \\
\hline
Capped $\ell_1$-penalty    & $\lambda\min(|t|, \theta)$, $(\theta > 0)$\\
\hline
\end{tabular*}
\end{sc}
\end{small}
\end{center}
\vskip -0.1in
\end{table}

\section{Majorization Minimization Algorithms}\label{sec:MM}

We consider a minimization problem with the objective function $F(\w)$. Given an estimate $\w^{(k)}$  at the $k$th iteration,
a typical MM algorithm consists of the following two steps:

\begin{itemize}
\item[] {\bf Majorization Step}:  Substitute $F(\w)$ by a tractable surrogate function $Q(\w|\w^{(k)})$,   such that
\[
Q(\w|\w^{(k)})\geq F(\w)
\]
for any $\w \in \dom F$, with equality  holding at $\w=\w^{(k)}$. 

\item[] {\bf Minimization Step}: Obtain the next parameter estimate $\w^{(k+1)}$  by minimizing $Q(\w|\w^{k})$ with respect to $\w$. That is,
\[
\w^{(k+1)}=\argmin_{\w} \; Q(\w | \w^{(k)}).
\] 
\end{itemize}


In order to address the global convergence of MM for solving the problem in (\ref{eq:problem}), we propose a generic MM framework under the assumptions given in the previous section.
We particularly present criteria to devise the majorant functions of the loss function and penalty function, respectively.


\subsection{Majorization of Loss Function}

We first consider the majorization of the loss function $f(\w)$.
Recall that  $\nabla f(\w)$  is assumed to be Lipschitz continuous
({Assumption} \ref{ASS1}).
Given the estimate $\w^{(k)}$ of $\w$ at the $k$th iteration, one would derive the majorization of $f(\w)$  to obtain $\w^{(k+1)}$.
For the sake of simplicity, we denote the corresponding surrogate function as $Q_f(\w|\w^{(k)})$.
In our work, we claim that $Q_f(\w| \w^{(k)})$ should
have the following two properties so that the surrogate can approximate the objective $f$ well and also lead to efficient computations.
\begin{assumption}\label{ASS:Mloss}
Let  $Q_f(\w|\w^{(k)})$ be the majorization  of $f(\w)$ such that
$Q_f(\w|\w^{(k)})\geq f(\w) $  and $Q_f(\w^{(k)}|\w^{(k)})=f(\w^{(k)})$. Additionally,
the following properties also hold:
\begin{enumerate}
\renewcommand{\labelenumi}{(\theenumi)}
\item[\emph{(i)}] $Q_f(\w|\w^{(k)})-f(\w)$ is $\gamma$-strongly convex, where $\gamma>0$;
\item[\emph{(ii)}] $\nabla Q_f(\w|\w^{(k)})$ is Lipschitz continuous, and $\nabla Q_f(\w^{(k)}|\w^{(k)})=\nabla f(\w^{(k)})$.
\end{enumerate}
\end{assumption}

Let us see several extant popular algorithms which meet Assumption~\ref{ASS:Mloss}.
Proximal algorithms \citep{rockafellar1976monotone,lemaire1989proximal,iusem1999augmented,combettes2011proximal,parikh2013proximal}
solve  optimization problems by using a so-called proximal operator of the objective function.
Suppose we have an objective function $f(\w)$ at hand.
Given the $k$th estimate $\w^{(k)}$, the proximal algorithm aims to solve the following problem
\[
\w^{(k+1)}=\argmin_{\w} \bigg\{f(\w)+\frac{1}{2 \alpha_k}\|\w-\w^{(k)}\|^2\bigg\},
\]
where $\alpha_k$ is the step size at each iteration. Typically,  $\w^{(k+1)}$ is written as:
\[
\w^{(k+1)}\denote \mathrm{Prox}_{\alpha_k f}(\w^{(k)}).
\]
The majorization function is defined as $Q_f(\w|\w^{(k)})\triangleq f(\w)+\frac{1}{2\alpha_k}\|\w-\w^{(k)}\|^2$.
When $f(\w)$ is convex and $\nabla f(\w)$ is Lipschitz continuous, it is easy to check that $Q_f(\w|\w^{(k)})$
satisfies Assumption~\ref{ASS:Mloss}.

Another powerful algorithm is the proximal gradient algorithm. The
algorithm is more efficient when dealing with the following problem
\begin{equation} \label{proxG}
\w^{*}=\argmin_{\w} \; \Big\{f(\w)+r(\w)\Big\},
\end{equation}
where $f$ is differentiable and convex and $r$ is nonsmooth. The proximal gradient algorithm first approximates $f(\w)$
based on a local linear expansion plus a proximal term, both at the current estimate $\w^{(k)}$. That is,
\[
f(\w)\approx f(\w^{(k)})+\langle \nabla f(\w^{(k)}), \w-\w^{(k)}\rangle + \frac{1}{2\alpha_{k}}\|\w-\w^{(k)}\|^2,
\]
where $\alpha_{k}>0$ is the step size.
Then the $(k{+}1)$th estimate of $\w$ is given as
\begin{equation}\label{eq:proxG}
\begin{split}
\w^{(k+1)}&=\argmin_\w \; \bigg\{ f(\w^{(k)})+\langle \nabla f(\w^{(k)}), \w-\w^{(k)}\rangle\\
& \qquad \qquad \quad + \frac{1}{2\alpha_{k}}\|\w-\w^{(k)}\|^2 + r(\w) \bigg\}.
\end{split}
\end{equation}
Equivalently,
\[
\w^{(k+1)}={\mathrm{Prox}}_{\alpha_k r}(\w^{(k)}-\alpha_k \nabla f(\w^{(k)})).
\]
Intuitively, the proximal gradient algorithm would take the gradient descent step first and then does the proximal minimization step.
In this algorithm, $Q_f(\w|\w^{(k)})\triangleq  f(\w^{(k)})+\langle \nabla f(\w^{(k)}), \w-\w^{(k)}\rangle + \frac{1}{2\alpha_{k}}\|\w-\w^{(k)}\|^2$. It is also immediately  verified that  $Q_f(\w|\w^{(k)})$
satisfies Assumption~\ref{ASS:Mloss} when $\alpha_k<L_{f}^{-1}$, where $L_{f}$ is the Lipschitz constant of $\nabla f(\w)$.

\begin{figure}[!ht]\label{fig:logist}
\vskip 0.2in
\begin{center}
\centerline{\includegraphics[width=80mm]{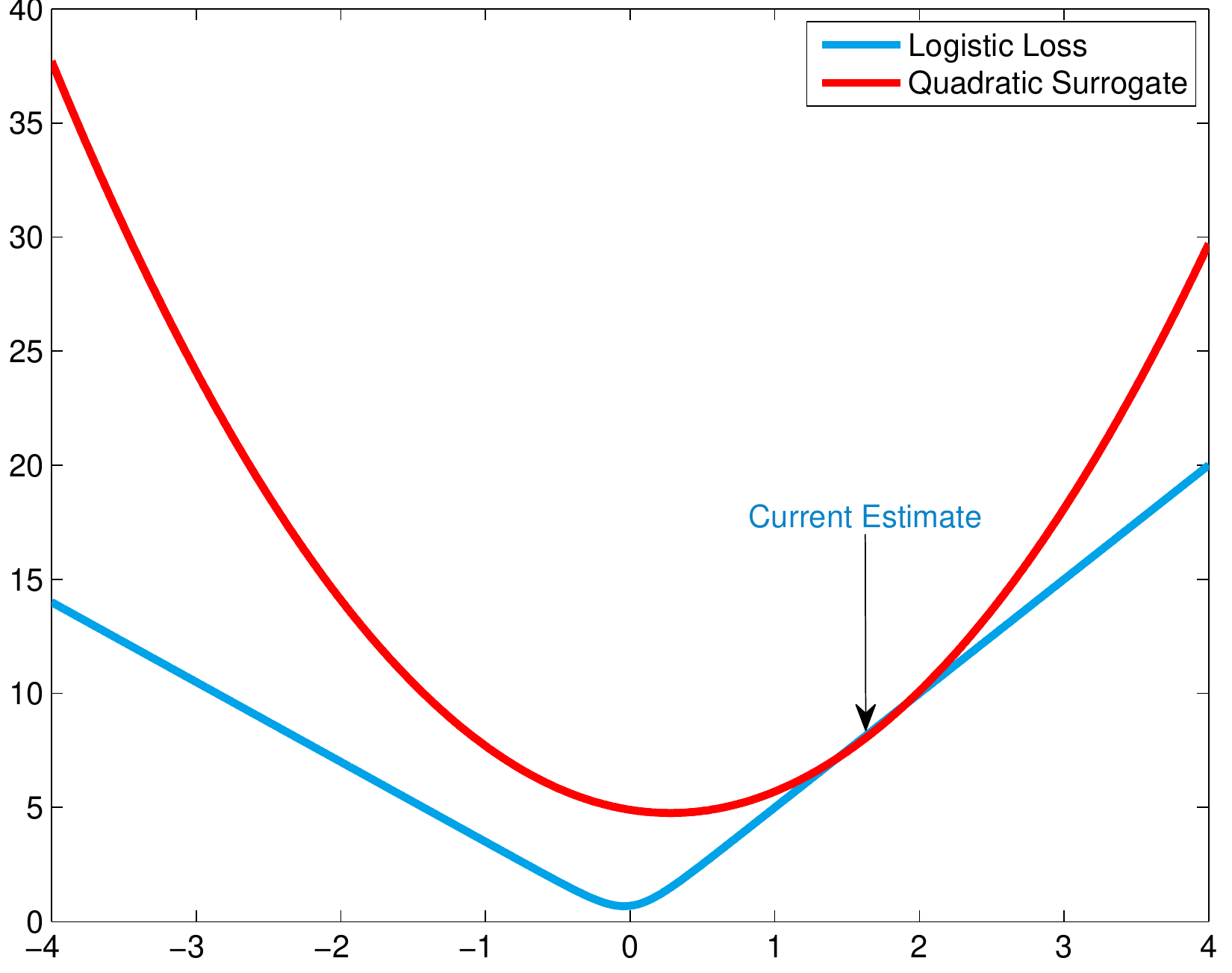}}
\caption{Surrogate for logistic loss}
\label{icml-historical}
\end{center}
\vskip -0.2in
\end{figure}

In fact, Lemma~\ref{lemm:1} implies that there always exists a quadratic surrogate of $f$ only if (\ref{A1}) holds. In particular,
we can define $Q_f$ as
\begin{equation} \label{Surr}
\begin{split}
Q_f(\w|\w^{(k)})=f(\w^{(k)})+\langle \nabla f(\w^{(k)}), \w-\w^{(k)}\rangle +\frac{\mu^{(k)}}{2}\|\w-\w^{(k)}\|^2,
\end{split}
\end{equation}
where we require that $\mu^{(k)}>L_{f}$.

\subsection{Majorization for Nonconvex Penalty Functions}

We assume that the penalty function $r(\w)=\sum_{i=1}^p \zeta(|w_i|)$.
Thus we can construct the surrogates for $\zeta$ separately.
It should be emphasized
that the majorization of the penalty function is not always necessary. For instance,
when one can easily obtain
\begin{equation}\label{Op1}
\w^{(k+1)}=\argmin_{\w} \bigg\{Q_f(\w|\w^{(k)})+r(\w)\bigg\},
\end{equation}
the surrogate for $r(\w)$ may not be considered. That is to say, this procedure is optional. However,  the surrogate for $r(\w)$ can result in efficient computations sometimes, especially when handling the proximal operator of $r(\w)$ suffers a large computation burden.

\begin{figure}[!ht]
\vskip 0.2in
\begin{center}
\centerline{\includegraphics[width=80mm]{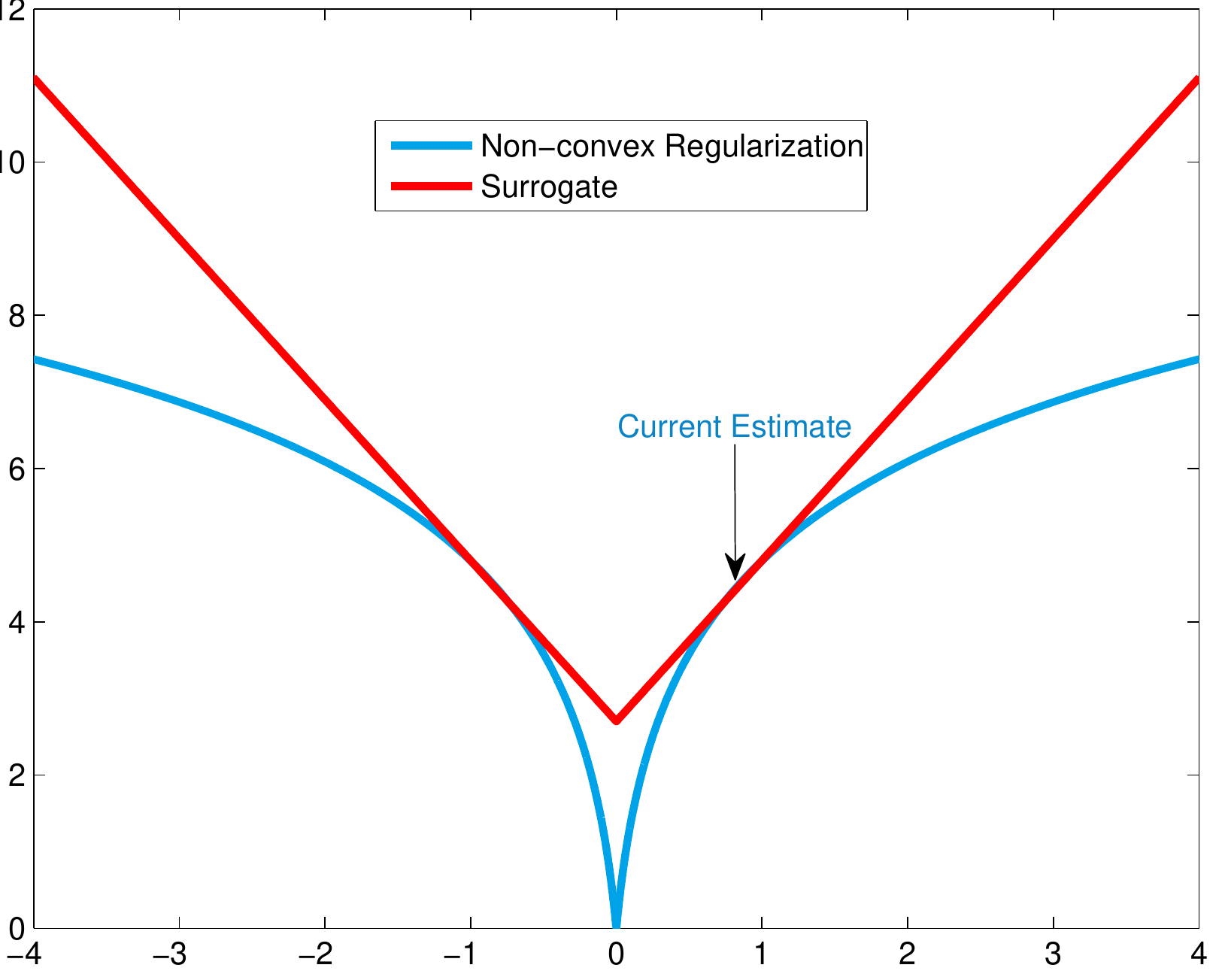}}
\caption{Surrogate for nonconvex regularization}
\label{fig:non-convex}
\end{center}
\vskip -0.2in
\end{figure}
We consider a more general case and  give some assumptions.
\begin{assumption}\label{ASS:Mreg}
Let $r(\w)=\sum_{i=1}^p \zeta(|w_i|)$, where the map $\zeta$: $\RB_+\to\RB_+$ is concave
and differentiable. Moreover, $\zeta^{\prime}(t)$ is Lipschitz continuous on $[0,+\infty)$. That is,
\[
|\zeta^{\prime}(t_1)-\zeta^{\prime}(t_2)|\leq L_{\zeta}|t_1-t_2|,
\]
for any $t_1, t_2 \geq 0$.
\end{assumption}

Many nonconvex penalties admit such properties, such as nonconvex LOG penalty, MCP, SCAD, etc. 
Although the $\ell_q$-norm ($q\in (0, 1)$) may not satisfy the gradient Lipschitz continuous condition, we alternatively consider $\zeta(|w|)=\lambda(1+\alpha|w|)^q$, with $\alpha>0$, which
is gradient Lipschitz continuous on $[0,+\infty)$.

Thanks to  concavity, we have
\begin{equation}\label{eq:Concave}
\zeta(|w_i|)\leq\zeta(|w^{(k)}_i|)+\zeta^{\prime}(|w^{(k)}_i|) (|w_i|-|w^{(k)}_i|),
\end{equation}
for any $i\in\{1,\cdots,p\}$. Thus, the majorant function for $r(\w)$, denoted by $Q_r(\w|\w^{(k)})$, is
\begin{equation}\label{eq:Mreg}
Q_r(\w|\w^{(k)})=\sum_{i=1}^{p}\Big[\zeta(|w^{(k)}_i|)+\zeta^{\prime}(|w^{(k)}_i|)(|w_i|-|w^{(k)}_i|)\Big].
\end{equation}
It is easy to see that $Q_r(\w|\w^{(k)})\geq r(\w)$ and $Q_r(\w^{(k)}|\w^{(k)})=r(\w^{(k)})$. Moreover, the corresponding surrogates transfer nonconvex objectives into convex ones, which brings  efficient and stable computations. As illustrated in Figure~\ref{fig:non-convex}, at each iteration, we optimize the tangent above the nonconvex penalty which is tight at the current estimate.

The key idea was early studied in DC programming \citep{gasso2009recovering}, which  linearizes iteratively concave functions to obtain convex surrogates. The idea has been also revisited by \cite{zou2008one,candes2008enhancing,chartrand2008iteratively}.

Specifically,
\cite{zou2008one} developed the local linear approximation (LLA) algorithm and pointed that the LLA algorithm can be
cast as an EM algorithm under certain condition. The LLA algorithm uses the same majorant function as in (\ref{eq:Mreg}) for the nonconvex and nonsmooth penalty function.

\cite{candes2008enhancing} studied  a so-called iteratively re-weighted $\ell_1$ minimization, which also falls into the MM procedure.
For example,
when  \[
r(\w)=\lambda\sum_{i=1}^p \log (1+\frac{1}{\epsilon}|w_i|), \mbox{ where }  \epsilon>0,
\]
the re-weighted $\ell_1$ minimization scheme is given as
\begin{equation}\label{eq:rel1}
  \w^{(k+1)}=\argmin_{\w} \bigg\{f(\w)+\lambda\sum_{i=1}^p\frac{|w_i|}{|w^{(k)}_i|+\epsilon}\bigg\},
\end{equation}
which can be also derived from (\ref{eq:Mreg}). 

 To be the best of our knowledge, there is few complete convergence results for these algorithms. In particular, it is hard to address the convergence of the sequence $\{\w^{(k)}\}_{k\in\mathbb N}$. 
In most of traditional treatments, the asymptotic stationary point is studied. These treatments usually follow the general convergence
results for MM algorithms \citep{lange2004optimization} without exploiting the property of the objective function.
Our global convergence results given in Section~\ref{sec:Convergence} are based on theory of  the Kurdyka-\L{}ojasiewicz inequality. Our results directly apply to the LLA and iteratively re-weighted $\ell_1$ minimization algorithms.

\subsection{A Generic MM Algorithm}
\begin{algorithm}[tb]
   \caption{Majorization Minimization Algorithm for Nonconvex Penalization}
   \label{alg:GNMM}
\begin{algorithmic}[1]
   \STATE Initialize $\w^{(0)}$ and $T$ (the maximum number of iterations), set $k=0$.
   \REPEAT

   \STATE Compute $Q_f(\w|\w^{(k)})$, which is the surrogate of $f(\w)$;
   \STATE Compute $Q_r(\w|\w^{(k)})$, which is the surrogate of $r(\w)$;
   \STATE Update $\w^{(k+1)}$ by (\ref{Op1}) or (\ref{Op2});
   \STATE $k = k + 1$;
   \UNTIL{Some stopping criterion is satisfied.}
\end{algorithmic}
\end{algorithm}

We are now ready to summarize the whole MM procedure.
Recall that the original problem (\ref{eq:problem}) includes two parts.
We first consider the simple case.
The ``simple" means the following problem can be handled easily:
\begin{equation}\label{eq:Simple}
\w^{*}=\argmin_\w \; \bigg\{\frac{1}{2\lambda}\|\w-\u\|^2+ r(\w) \bigg\}.
\end{equation}
This implies that (\ref{Op1}) can be efficiently solved. This leads to nonconvex proximal-gradient methods \citep{fukushima1981generalized,lewis2008proximal}.
We thus generate a sequence $\{\w^{(k)}\}_{k\in \NB}$ by (\ref{Op1}).

However, when (\ref{eq:Simple}) is intractable, we substitute $r(\w)$ with (\ref{eq:Mreg}).
Then in each iteration, the problem reduces to
\begin{equation}\label{Op2}
\w^{(k+1)}=\argmin_\w \;  \bigg\{Q_f(\w|\w^{(k)})+Q_r(\w|\w^{(k)}) \bigg\}.
\end{equation}
The whole procedure is briefly presented in Algorithm \ref{alg:GNMM}.

\section{Convergence Analysis}\label{sec:Convergence}


We  now study the convergence analysis of  Algorithm \ref{alg:GNMM}.
It should be claimed that the global convergence, which is our focus, means that for any $\w^{(0)}\in\RB^p$, the sequence $\{\w^{(k)}\}_{k\in\mathbb N}$ generated by (\ref{Op1}) or (\ref{Op2}) converges to the critical point of $F(\w)$.

\begin{lemma}\label{lem:descent}
Suppose {\bf Assumptions  \ref{ASS3}, \ref{ASS:Mloss}} hold or {\bf Assumptions \ref{ASS3}, \ref{ASS:Mloss}, and \ref{ASS:Mreg}} hold. Then the sequence
$\{\w^{(k)}\}_{k\in\mathbb N}$  generated by (\ref*{Op1}) or generated by (\ref*{Op2}) satisfies the following properties.
\begin{enumerate}
\item[\emph{(i)}] The generated sequence $\left\{F(\w^{(k)})\right\}_{k\in \mathbb N}$ is non-increasing, specifically,
\[F(\w^{(k)})-F(\w^{(k+1)})\geq \frac{\gamma}{2}\|\w^{(k)}-\w^{(k+1)}\|^2, \quad \forall k \geq 0.
    \]
\item[\emph{(ii)}] \[\sum_{k=0}^{\infty}\|\w^{(k+1)}-\w^{(k)}\|^2<+\infty, \]
which implies $\lim_{k\to\infty}\|\w^{(k+1)}-\w^{(k)}\|=0$.
\end{enumerate}
\end{lemma}


Lemma \ref{lem:descent}  enjoys the descent property of the MM approach which always makes the objective term decrease after each iteration. Moreover, the objective function value decreases at least $\frac{\gamma}{2}\|\w^{(k+1)}-\w^{(k)}\|^2$ for the $k$th step.
By the fact that $\inf_{\w} F(\w)>-\infty$, we can draw the conclusion that the sequence $\{F(\w^{(k)})\}_{k\in\NB}$ converges.

Because  of the coerciveness of function $F(\w)$, there exists a convergent subsequence $\left\{\w^{(n_k)}\right\}_{k\in\mathbb N}$ that converges to $\bar\w$.
The set of all cluster or limit points which are started with $\w^{(0)}$ is denoted by ${\mathcal M}(\w^{(0)})$. That is,
\[
{\mathcal M}(\w^{(0)})\denote \bigg\{\bar \w \in \RB^p: \exists  ~n_k,~\left\{n_k\right\}_{k \in \NB}, \mbox{~such that }~\w^{(n_k)}\to\bar \w \mbox{~as} ~k ~\to ~\infty\bigg\}.
\]

It is also easy to see that $F(\w)$ is constant and finite on $\mathcal M(\w^{(0)})$. In the following lemma we attempt to demonstrate that all points which belong to $\mathcal M(\w^{(0)})$ are stationary or critical points of $F(\w)$.

\begin{lemma}\label{lem:bound-sub}
Suppose {\bf Assumptions \ref{ASS1}, \ref{ASS3}, \ref{ASS:Mloss}} hold, and the sequence
$\{\w^{(k)}\}_{k \in \NB}$ is generated by (\ref*{Op1}). Let
$A^{(k+1)}=\nabla f(\w^{(k+1)})-\nabla Q_f(\w^{(k+1)}| \w^{(k)})$. Then
\begin{enumerate}
\item[\emph{(i)}]  $A^{(k+1)}\in \partial F(\w^{(k+1)})$;
\item[\emph{(ii)}]  $\|A^{(k+1)}\|\leq (L_{Q_f}+L_{f})\|\w^{(k+1)}-\w^{(k)}\|$.
\end{enumerate}
\end{lemma}

\begin{lemma}\label{lem:share}
Suppose $r(\w)=\sum^p_{i=1}\zeta(w_i)$, where $\zeta$: $\RB_+\to\RB_+$ is concave
and  continuous differentiable on $[0,+\infty)$. Let $Q_r(\w|\w^{(k)})=\sum_{i=1}^{p}\zeta(|w^{(k)}_i|) + \zeta^{\prime}(|w^{(k)}_i|)(|w_i|-|w^{(k)}_i|)$. Then
\begin{enumerate}
\item[\emph{(i)}]  $Q_r(\w|\w^{(k)})\geq r(\w)$ \mbox{and} $Q_r(\w^{(k)}|\w^{(k)})=r(\w^{(k)})$;
\item[\emph{(ii)}]  $\partial Q_r(\w^{(k)}|\w^{(k)})=\partial r(\w^{(k)})$.
\end{enumerate}
\end{lemma}

Lemma \ref{lem:share} shows the relationship between the nonconvex (nonsmooth) penalty function and the corresponding surrogate.  This also implies that the surrogate approximates the penalty function well.

We introduce the notion of  $\sgn(u)$, which is defined as

\begin{equation} \label{def: sgn}
\sgn(u) \denote \left\{
\begin{array}{lc}
1 & \mbox{ if } u > 0, \\
c & \mbox{ if } u = 0, \\
-1 & \mbox{ if } u < 0.
\end{array}
\right.
\end{equation}
Here $c$ is some  real number in  $[-1, 1]$.  
We emphasize that $\sgn(u)$ is a scalar rather than a set.

\begin{lemma}[Main Lemma]\label{cor:bound-sub}
Suppose {\bf Assumptions \ref{ASS1}, \ref{ASS3}, \ref{ASS:Mloss}, \ref{ASS:Mreg}} hold, and the sequence
$\{\w^{(k)}\}_{k\in \NB}$ is generated by (\ref*{Op2}). Let $b^{(k)}_i= \sgn(w_i^{(k+1)})[\zeta^{\prime}(|w_i^{(k)}|)-\zeta^{\prime}(|w_i^{(k+1)}|)]$ for $i\in[p]$, and
$B^{(k+1)}=\nabla f(\w^{(k+1)})-\nabla Q_f(\w^{(k+1)}| \w^{(k)})- (b^{(k)}_1,b^{(k)}_2, \cdots, b^{(k)}_p)^T$. Then
\begin{enumerate}
\item[\emph{(i)}]  $B^{(k+1)}\in \partial F(\w^{(k+1)})$;
\item[\emph{(ii)}]  $\|B^{(k+1)}\|\leq (L_{Q_f}+L_{f}+L_{\zeta})\|\w^{(k+1)}-\w^{(k)}\|$.
\end{enumerate}
\end{lemma}

Both Lemma \ref{lem:bound-sub} and Lemma \ref{cor:bound-sub} suggest a subgradient lower bound
for the iterate gap. Due to the majorization of the nonconvex and nonsmooth penalty functions, it is more challenging to bound the subgradient. The ingredient is to observe that the majorant function and the original one share the same subgradient at the current estimate. With Lemmas \ref{lem:bound-sub} and \ref{cor:bound-sub}, we are now ready to state the following lemma.

\begin{lemma}\label{lemm:crictical}
Suppose {\bf Assumptions \ref{ASS1}, \ref{ASS3}, \ref{ASS:Mloss}, \ref{ASS:Mreg}} hold.
Let the sequence $\left\{\w^{(k)}\right\}_{k \in \NB}$ be generated by (\ref{Op1})
 or (\ref{Op2}). Then
\begin{enumerate}
\item[\emph{(i)}] $\mathcal M(\w^{(0)})$ is not empty and ${\mathcal M}(\w^{(0)}) \subset \mathrm{crit} F$;
\item[\emph{(ii)}] \begin{equation}\label{eq:set}
\lim_{k \to \infty} \mathrm{dist} \Big(\w^{(k)}, {\mathcal M}(\w^{(0)}\Big)=0.
\end{equation}
\end{enumerate}
\end{lemma}

Lemma \ref{lemm:crictical} implies that   ${\mathcal M}(\w^{(0)})$ is the subset of stationary or critical points of $F(\w)$ and $\{\w^{(k)}\}_{k\in\NB}$  are approaching to one point of ${\mathcal M}(\w^{(0)})$.   Our current concern is to prove $\lim_{k\to\infty}\w^{(k)}=\w^{*}$. From  \cite{lange2004optimization}, we know that $\mathcal M(\w^{(0)})$ is connected. Additionally, if $\mathcal M(\w^{(0)})$ is finite, $\{\w^{(k)}\}_{k\in \NB}$ converges .

We can obtain the global convergence based on the assumption that $\mathcal M(\w^{(0)})$ is finite. However, the assumption is not practical, because it is usually unknown. Moreover, it is hard  to  check this assumption. To avoid this issue, the Kurdyka-\L{}ojasiewicz property of the objective function enters in action, because it is often a very easy task to verify the Kurdyka-\L{}ojasiewicz property of a function.

\begin{theorem}\label{GlobalCon}
Suppose that $F$ has Kurdyka-\L{}ojasiewicz property at each point of $\dom \, \partial F$, and {\bf Assumptions \ref{ASS1}, \ref{ASS3}, \ref{ASS:Mloss}, \ref{ASS:Mreg}} hold. Let the sequence $\left\{\w^{(k)}\right\}_{k\in \mathbb N}$ be generated by scheme (\ref{Op1}) or (\ref{Op2}). Then the following assertions hold.
\begin{enumerate}
\item[\emph{(i)}] The sequence $\left\{\w^{(k)}\right\}_{k \in \NB}$ has finite length.
\begin{equation}\label{eq:gloabl}
\sum^\infty_{k=0}\Big\|\w^{(k+1)}-\w^{(k)}\Big\|<\infty
\end{equation}
\item[\emph{(ii)}] The sequence $\left\{\w^{(k)}\right\}_{k \in \NB}$ converges to a critical point $\w^*$of F.
\end{enumerate}
\end{theorem}

Theorem \ref{GlobalCon} shows the global convergence of Algorithm \ref{alg:GNMM}. As we have stated,
many methods for solving a nonconvex  and nonsmooth problem, such as the re-weighted $\ell_1$ \citep{candes2008enhancing}
and LLA \citep{zou2008one},  share the same convergence property
as in Theorem \ref{GlobalCon}.
\citet{attouch2010proximal,bolte2013proximal}  have  well established  the global convergence for nonconvex and nonsmooth problems based on the Kurdyka-\L{}ojasiewicz  inequality. It is also interesting to point out that their procedures
fall into (\ref{Op1}). However, they focused on a coordinate descent procedure. The work  of
\citet{attouch2010proximal,bolte2013proximal} cannot be trivially extended to our more general case.

\begin{theorem}[Convergence Rate]\label{theorem:convergenceRate}
Suppose {\bf Assumptions \ref{ASS1}, \ref{ASS3}, \ref{ASS:Mloss}, \ref{ASS:Mreg}} hold, and $\{\w^{(k)}\}_{k\in\mathbb N}$ is generated
by (\ref{Op1}) or (\ref{Op2})  which converges to a critical point $\w^*$ of $F$, which satisfies the Kurdyka-\L{}ojasiewicz  property at each point of $\dom \, \partial F$ with $\phi(t)=ct^{1-\theta}$ for $c>0$ and $\theta\in[0,1)$. We have
\begin{enumerate}
\item[\emph{(i)}] if $\theta=0$, $\{\w^{(k)}\}_{k\in\mathbb N}$ converges to $\w^*$ in finite iterations;
\item[\emph{(ii)}] if $\theta\in(0,\frac{1}{2}]$, $\|\w^{(k)}-\w^*\|\leq C\rho^k$, $\forall k\geq K_{0}$, for some $K_0>0, C >0,\rho\in(0,1)$;
\item[\emph{(iii)}] if $\theta\in(\frac{1}{2},1)$, $\|\w^{(k)}-\w^*\|\leq Ck^{-\frac{1-\theta}{2\theta-1}}$,$\forall k\geq K_{0}$, for some $K_0>0, C >0$.
\end{enumerate}
\end{theorem}

Theorem \ref{theorem:convergenceRate} tells us the convergence rate of our MM procedure for solving the nonconvex regularized problem, which is based on the geometrical property  of the function $F$ around its critical point. We see that the convergence rate is at least sublinear.




\section{Extension to Concave-Convex Procedure}\label{sec:cccp}


In this section we show that our work can be extended to the concave-convex procedure (CCCP) \citep{yuille2003concave}.
It is worth noting that CCCP can be also unified into the MM framework.  

The CCCP is usually used to solve the following problem:
\begin{equation}\label{eq:CCCP}
\begin{split}
\min_{\w} ~~~~~~~~~~&u(\w)-v(\w), \\
\st~~~~~~~~~~~ &c_{i}(\w)\leq\0, ~i \in[n], \\
&d_{j}(\w)=\0, ~j\in[m],
\end{split}
\end{equation}
where $u$, $v$ and $c_i$ are real-valued convex functions and $d_j$ are affine functions.
The CCCP algorithm aims to solve the following sequence of convex optimization problems:
\begin{equation}\label{eq:Convex-CCCP}
\begin{split}
\w^{(k+1)} = \argmin_{\w} ~~~~~~~ &u(\w)-\nabla v(\w^{(k)})^T \w, \\
\st ~~~~~~~~ &c_{i}(\w)\leq\0, ~i\in[n], \\
&d_{j}(\w)=\0, ~j\in[m].
\end{split}
\end{equation}
Denote $\mathcal C=\Big\{\w:c_{i}(\w)\leq\0, d_{j}(\w)=\0, i\in[n], j\in[m]\Big\}$, and let
$\delta_{\mathcal C}(\w)$ be the indicator function of the feasible set $\mathcal C$; that is,
\[
\delta_{\mathcal C}(\w)=\left\{
\begin{array}{c}
0,~~~~~~~~~~~~~~~~~\w\in\mathcal C, \\
+\infty,~~~~~~~~~~~~~\w \notin \mathcal C.
\end{array}
\right.
\]
It is directly proved that $\delta_{\mathcal C}$ is a convex function.
Now the original problem can be reformulated as
\begin{equation}\label{eq:Indicator_cccp}
\min_{\w} \;  F(\w) \triangleq \Big\{\delta_{\mathcal C}(\w)+u(\w)-v(\w)\Big\}.
\end{equation}
Thus,  the CCCP approach would solve the following convex problem at each iteration:
\begin{equation}\label{eq:Indicator_ccp_surrogate}
\w^{(k+1)} = \; \argmin_{\w}~~~~\Big\{ \delta_{\mathcal C}(\w)+u(\w)-\nabla v(\w^{(k)})^T\w \Big\}.
\end{equation}

In fact, the CCCP approach can be viewed as an MM algorithm. In particular,
since $v(\w)$ is convex, $-v(\w)$ is concave.
As a result, we have
\[
-v(\w)\leq-v(\w^{(k)})-\nabla v(\w^{(k)})^T(\w-\w^{(k)}).
\]
This leads us to  the linear majorization  of $-v(\w)$. When the constant part is omitted,
(\ref{eq:Convex-CCCP}) or (\ref{eq:Indicator_ccp_surrogate}) are recovered. In summary,
CCCP  linearizes the concave part of the objective function.
Next, we make some assumptions to address the convergence of CCCP.
\begin{assumption}\label{ASS:CCCP}
Consider the problem in (\ref{eq:Indicator_cccp}) where $\delta_{\mathcal C}$, $u$,  and $v$
are  convex functions. Suppose the following three asserts  hold.
\begin{enumerate}
\renewcommand{\labelenumi}{(\theenumi)}
\item[\emph{(i)}]  $u(\w)$ and $v(\w)$ are $C^1$ functions;
\item[\emph{(ii)}]  u(\w) is $\gamma$-strongly convex;
\item[\emph{(iii)}]  $\nabla v(\w)$ is Lipschitz continuous.
\end{enumerate}
\end{assumption}

With the above assumption, the following theorem shows that
the sequence $\{\w^{(k)}\}_{k\in\mathbb N}$ generated by CCCP  converges to the critical point of $F(\w)$.

\begin{theorem}[Global Convergence of CCCP]\label{the:cccp}
Suppose {\bf Assumption \ref{ASS3}, \ref{ASS:CCCP}} hold.  And F satisfy  the
 the Kurdyka-\L{}ojasiewicz  property at each point of $\dom \, \partial F$. Let the sequence $\left\{\w^{(k)}\right\}_{k\in \mathbb N}$ be generated by  (\ref{eq:Indicator_ccp_surrogate}).  Then the conclusions of Theorems~\ref{GlobalCon}
and \ref{theorem:convergenceRate} hold.
\end{theorem}

It is worth pointing out that the global convergence analysis for CCCP has also been studied by \citet{lanckriet2009convergence}.
Their analysis is based on the novel Zangwill's theory. Zangwill's theory is  a very important
tool to deal with the convergence  issue of iterative algorithms. But it typically requires that $\mathcal M(\w^{(0)})$ is finite or discrete to achieve the convergent sequence $\{\w^{(k)}\}_{k\in\mathbb N}$ \citep{wu1983convergence,lanckriet2009convergence}.
In contrast,  our analysis based on the Kurdyka-\L{}ojasiewicz inequality does not need this requirement.

\section{Numerical Analysis}\label{sec:experiment}

In this paper our principal focus has been to explore the convergence properties of  majorization minimization (MM) algorithms
for nonconvex optimization problems. However, we have also developed two special MM algorithms based on (\ref{Op1}) and  (\ref{Op2}), respectively.
Thus, it is
interesting to conduct empirical analysis of convergence of the  algorithms. We particularly  employ the logistic loss and  LOG penalty for the classification problem.
We refer to the algorithms as MM-(a) and  MM-(b) for discussion simplicity.


We evaluate both MM-(a) and MM-(b) on binary datasets\footnote{http://www.csie.ntu.edu.tw/~cjlin/libsvmtools/datasets/binary.html}. Descriptions of the datasets are reported in Table \ref{tab1}.
For each dataset $\big\{\x_i,y_i\big\}_{i=1}^n$,
\begin{equation}\label{exp:obj}
F(\w)=\underbrace{\frac{1}{n}\sum_{i=1}^n\log(1+\exp(-y_i\x_i^T\w))}_{f(\w)}+\underbrace{\lambda\sum_{i=1}^p \log(1+\alpha|w_i|)}_{r(\w)},
\end{equation}
where $\lambda>0$ and $\alpha>0$ are hyperparameters.
We adopt the corresponding majorizantion function
\[
Q_f(\w|\w^{(k)})=f(\w^{(k)})+\langle \nabla f(\w^{(k)}), \w-\w^{(k)}\rangle +\frac{\mu^{(k)}}{2}\|\w-\w^{(k)}\|^2
\] and
\[ Q_r(\w|\w^{(k)})=\lambda \sum_{i=1}^p\bigg\{\log(1+\alpha |w^{(k)}_i|)+ \frac{\alpha}{1+\alpha |w^{(k)}_i|} (|w_{i}| - |w^{(k)}_i|) \bigg\}.
\]
As mentioned in the previous section,   the Lipschitz constant of $\nabla f(\w)$ is bounded by $\frac{1}{4n}\sum_{i=1}^n \x^T_{i}\x_{i}$. Typically, to set the value $\mu^{(k)}$, one often uses
the line-search method \citep{beck2009fast} to achieve better performance. However, since we are only concerned with  the convergence behavior of MM,
we just set $\mu^{(k)}= \frac{\rho}{4n}\sum_{i=1}^n \x^T_{i}\x_{i}$ where $\rho \geq 1$.

\begin{figure*}[!ht]
\begin{center}
\centering
\subfigure[leukemia $(\lambda=0.3,\alpha=0.3)$]{\includegraphics[width=70mm]{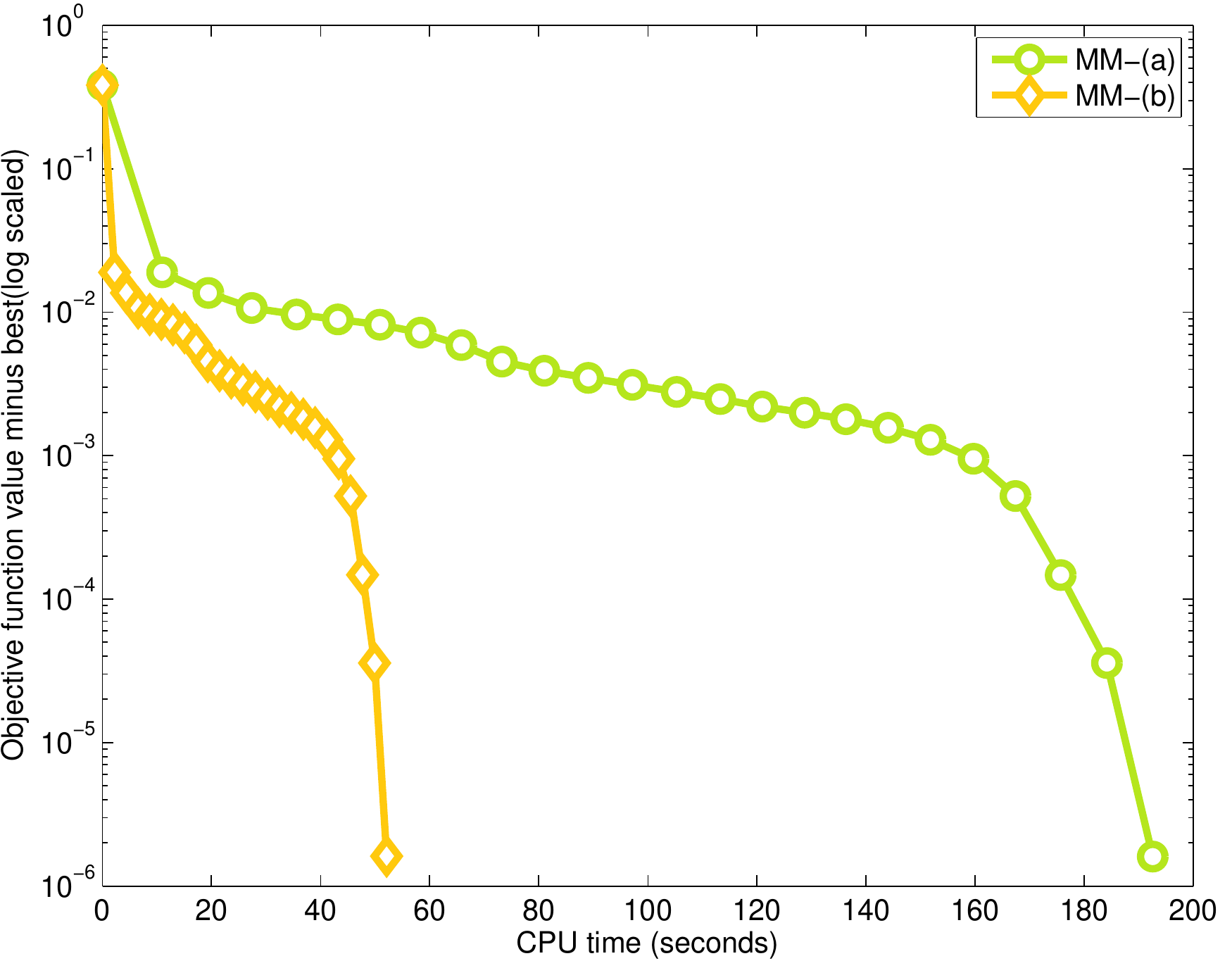}}
\subfigure[leukemia $(\lambda=0.0001, \alpha=0.0003)$]{\includegraphics[width=70mm]{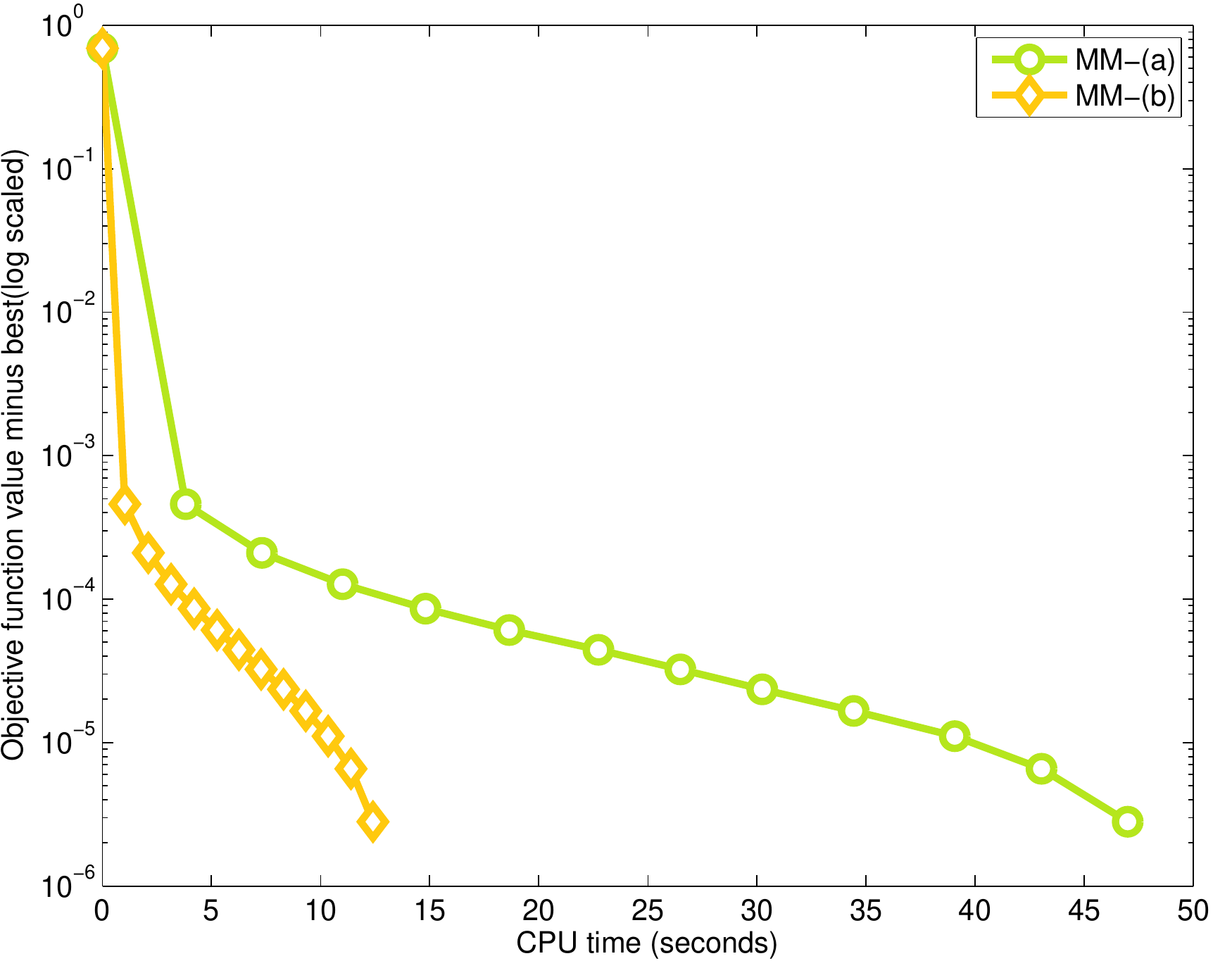}}
\subfigure[news20 $(\lambda=0.01,\alpha=0.03)$]{\includegraphics[width=70mm]{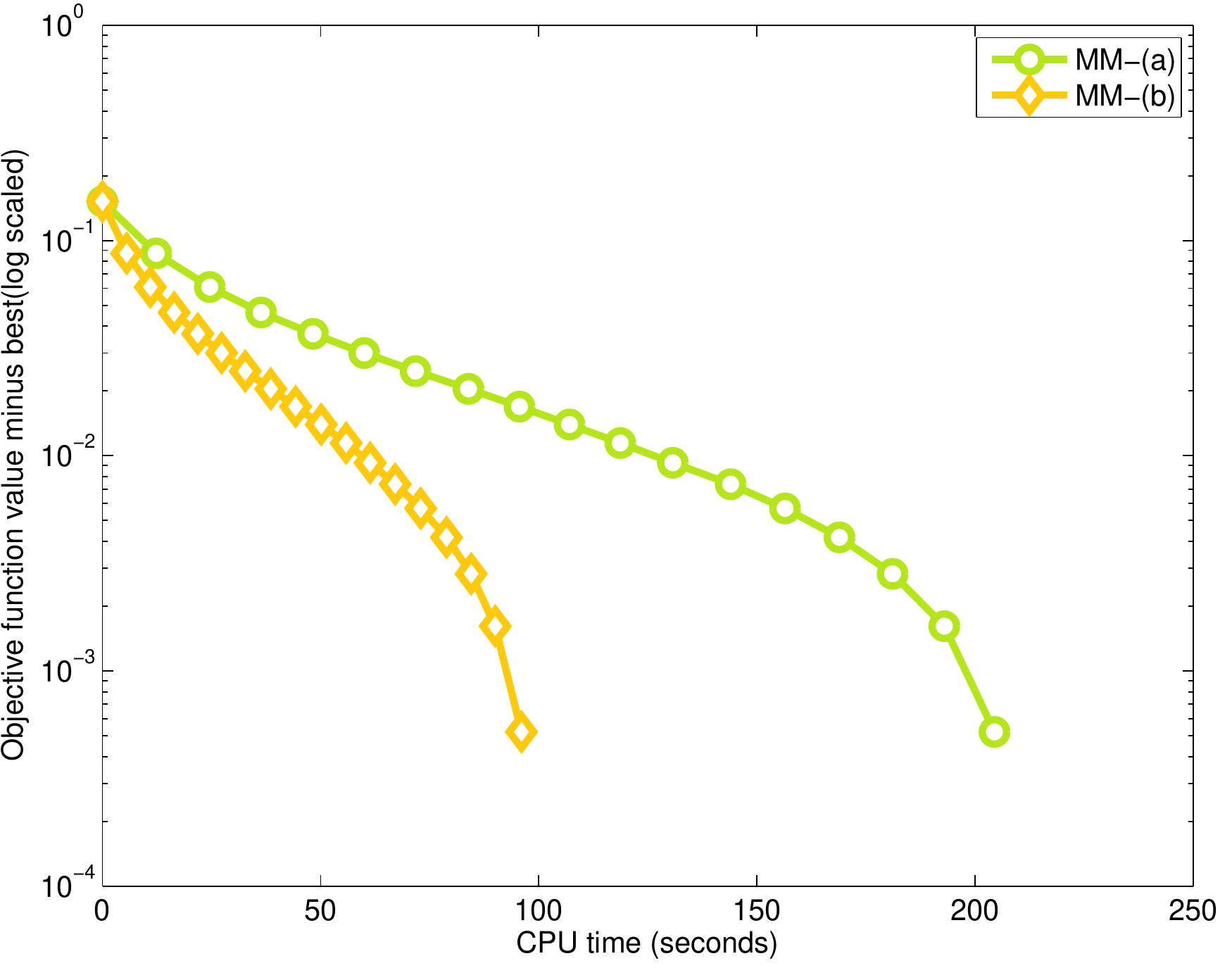}}
\subfigure[news20 $(\lambda=0.0001, \alpha=0.0003)$]{\includegraphics[width=70mm]{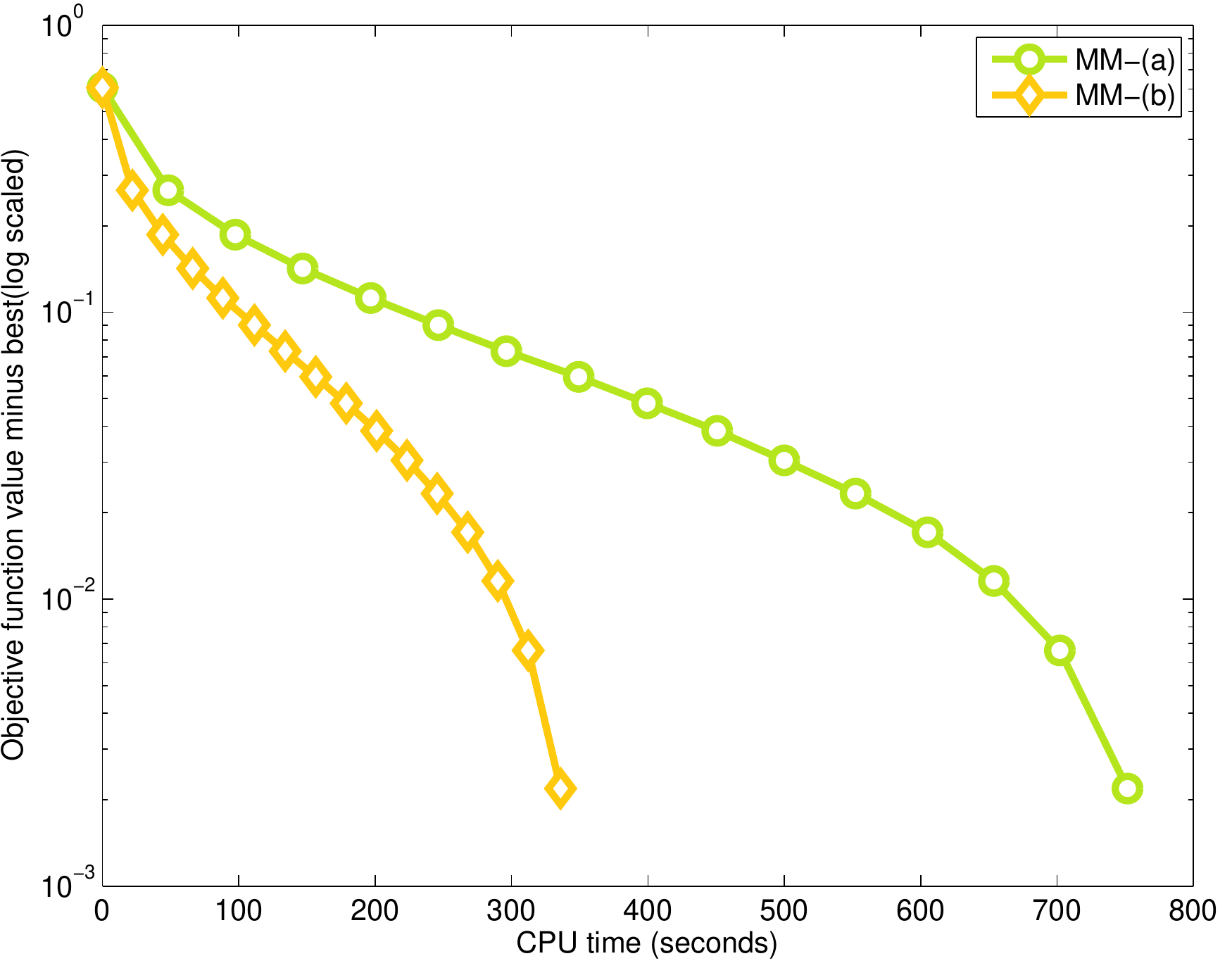}}
\subfigure[covtype $(\lambda=0.01, \alpha=0.03)$]{\includegraphics[width=70mm]{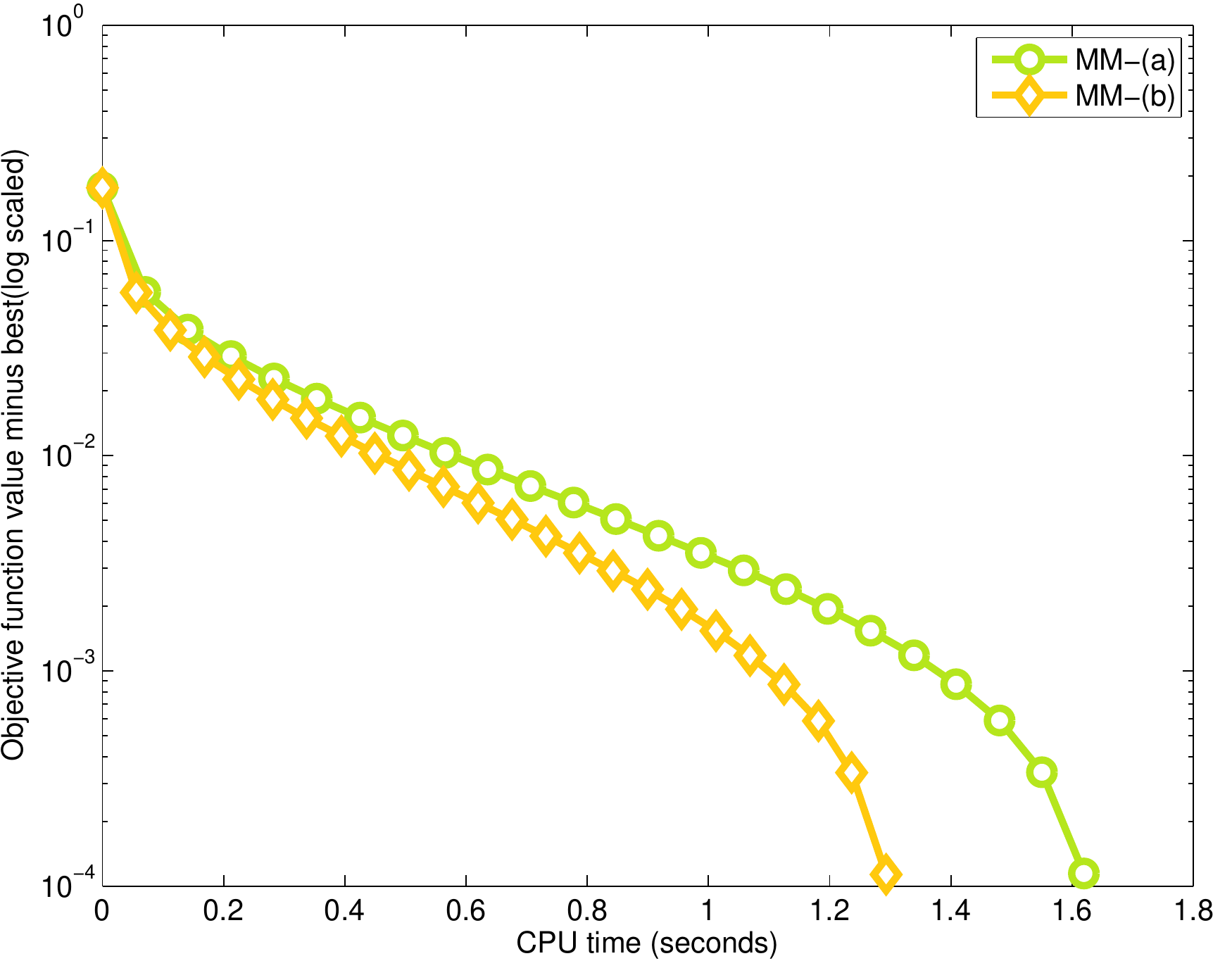}}
\subfigure[covtype $(\lambda=0.0001, \alpha=0.0003)$]{\includegraphics[width=70mm]{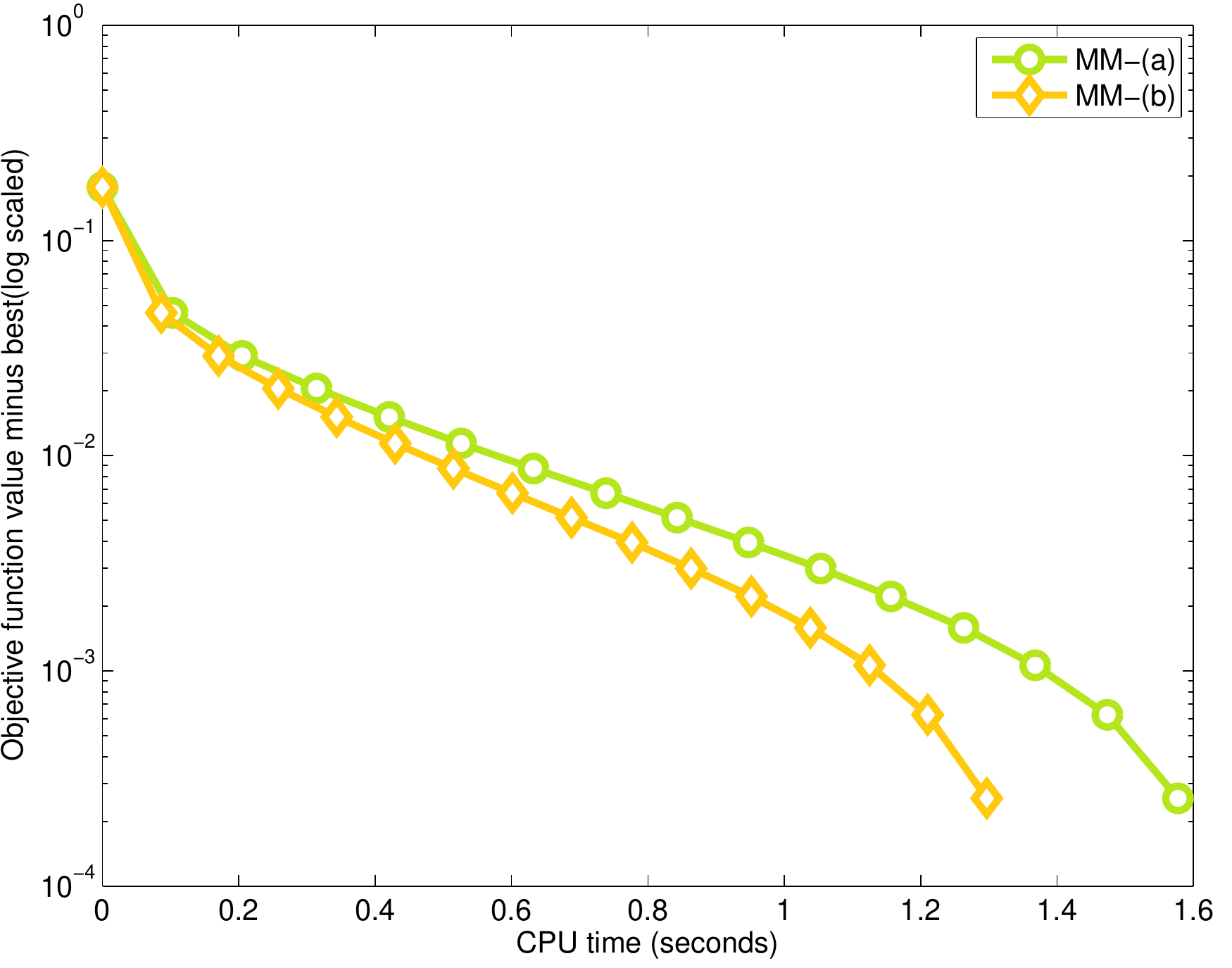}}
\end{center}
   \caption{performance of MM with different parameter settings}
\label{fig:obj2time}
\end{figure*}

We  plot the error between objective function values and the $F(\w^*)$($\log$ scaled) vs. CPU times with respect to different  hyperparameters settings
in Figure \ref{fig:obj2time}. We observe that both MM-(a) and MM-(b) generate the monotone decreasing sequence $\{F(\w^{(k)})\}_{k\in \mathbb N}$ and achieve nearly the same
optimal objective value. We also find that MM-(b) runs faster than MM-(a).
This implies that it is efficient to construct the majorization function of the LOG
penalty. In fact, MM-(a) will cost more computations when one directly calculates the proximal operator of the LOG penalty.
In contrast, MM-(b) only needs to do the soft-thresholding (shrinkage) operator on the current estimate. In summary, numerical experiments show that both
MM-(a) and MM-(b) make the objective function value decrease and converge.

\begin{table}[!ht]\setlength{\tabcolsep}{10pt}
\caption{Description of the datasets}
\label{tab1}
\begin{center}
\begin{footnotesize}
\begin{tabular}{|l|c|c|c|}
\hline
	{Data sets} & n& p & storage      \\
\hline
    {leukemia}   &  $72 $&$7129$&  sparse\\
\hline
    { news20}&   $19996  $   & $1355191$   &  sparse \\
\hline
    { covtype}&   $ 581012$  & $54$ &   dense\\
\hline

\end{tabular}
\end{footnotesize}
\end{center}
\end{table}


\section{Conclusions}\label{sec:conclusion}

Majorization minimization (MM) algorithms are very popular in machine learning and statistical inference. 
In this paper, we have employed MM algorithms to solve the nonconvex regularized problems. 
However, the convergence analysis of MM for nonconvex and nonsmooth problems is a challenging issue. 
We have established the global convergence results of the MM procedure using the geometrical property of the objective function. 
In particular, our results are built on the  Kurdyka-\L{}ojasiewicz inequality. 
We have shown that our results also apply to the iteratively re-weighted $\ell_1$ minimization method,  local linear approximation (LLA), 
and concave-convex procedure (CCCP).

\appendix

\section{Proofs} \label{sec:proof}

\subsection{The proof of Corollary \ref{cor:1}}
\begin{proof}
Since $f_i(\w)$ is differentiable for $i\in[n]$ and each $\nabla f_{i}(\w)$ is $L_i$-Lipschitz continuous, we have
\[
 \|\nabla f_i(\u)-\nabla f_i(\v)\|\leq L_{i}\|\u-\v\|, 
\]
for $i \in [n]$.
Then\[
\begin{split}
\|\nabla h(\u)-\nabla h(\v)\|&=\|\sum_{i=1}^{n}\alpha_i \nabla f_{i}(\u)-\sum_{i=1}^{n}\alpha_i \nabla f_{i}(\v)\| \\
&\leq\sum_{i=1}^n |\alpha_i|\|\nabla f_{i}(\u)-\nabla f_{i}(\v)\|\\
&\leq(\sum_{i=1}^n |\alpha_{i}|L_i)\|\u-\v\|
\end{split}
\]
So, $\nabla h(\w)$ is $\sum_{i=1}^n |\alpha_{i}|L_i$ Lipschitz continuous.
\end{proof}

\subsection{The proof of Lemma \ref{lem:descent}}
\begin{proof}
We first consider (\ref{Op1}) procedure.
\begin{enumerate}
\item[{(i)}]
Recall that 
\[
\w^{(k+1)}=\argmin_{\w} \bigg\{Q_f(\w|\w^{(k)})+r(\w)\bigg\}.
\]
We have
\begin{equation}\label{eq:Opt1Condition}
Q_f(\w^{(k+1)}|\w^{(k)})+r(\w^{(k+1)})-Q_f(\w^{(k)}|\w^{(k)})-r(\w^{(k)})\leq 0.
\end{equation}
By the strongly-convex property of $Q_f(\w|\w^{(k)})-f(\w)$,
\[
\begin{split}
Q_f(\w^{(k+1)}|\w^{(k)})&-f(\w^{(k+1)})- Q_f(\w^{(k)}|\w^{(k)})+f(\w^{(k)})\geq\\\langle \nabla Q_f(\w^{(k)}|\w^{(k)})-\nabla f(\w^{(k)}),&\w^{(k+1)}-\w^{(k)}\rangle+\frac{\gamma}{2}\|\w^{(k+1)}-\w^{(k)}\|^2\\
&=\frac{\gamma}{2}\|\w^{(k+1)}-\w^{(k)}\|^2
\end{split}
\]
The last equality complies with $\nabla Q_f(\w^{(k)}|\w^{(k)})=\nabla f(\w^{(k)})$.

Combining with (\ref{eq:Opt1Condition}), we have
\[
f(\w^{(k)})+r(\w^{(k)})-f(\w^{(k+1)})-r(\w^{(k+1)})\geq \frac{\gamma}{2}\|\w^{(k+1)}-\w^{(k)}\|^2.
\]
So
\[
F(\w^{(k)})-F(\w^{(k+1)})\geq \frac{\gamma}{2}\|\w^{(k+1)}-\w^{(k)}\|^2.
\]
\item[{(ii)}] We summary the above inequality from $k=0$ to $+\infty$. Then
\[
\sum_{k=0}^{+\infty} F(\w^{(k)})- F(\w^{(k+1)}) \geq \sum_{k=0}^{+\infty} \frac{\gamma}{2}\|\w^{(k+1)}-\w^{(k)}\|^2
\]
Notice $\inf_{\w} F(\w)>-\infty$, so 
\[
 \sum_{k=0}^{+\infty}\|\w^{(k+1)}-\w^{(k)}\|^2\leq\frac{2}{\gamma} (F(\w^{(0)})-F(\w^{(\infty)}))<+\infty,
\]
which completes the proof.
\end{enumerate}

Let's come to the sequence $\{\w^{(k)}\}_{k\in\mathbb N}$ generated by (\ref{Op2}).

Similarly,
\[
\w^{(k+1)}=\argmin_{\w} \bigg\{Q_f(\w|\w^{(k)})+Q_r(\w|\w^{(k)})\bigg\}.
\]
We obtain
\begin{equation}\label{eq:Opt2Condition}
Q_f(\w^{(k+1)}|\w^{(k)})+Q_r(\w^{(k+1)}|\w^{(k)})-Q_f(\w^{(k)}|\w^{(k)})-Q_r(\w^{(k)}|\w^{(k)})\leq 0.
\end{equation}
Since we have(similar to proof of Lemma \ref{lem:descent})
\[
Q_f(\w^{(k+1)}|\w^{(k)})-f(\w^{(k+1)})- Q_f(\w^{(k)}|\w^{(k)})+f(\w^{(k)})\geq \frac{\gamma}{2}\|\w^{(k+1)}-\w^{(k)}\|^2.
\]
On the other hand,
\[
\begin{split}
Q_r(\w^{(k+1)}|\w^{(k)})-Q_r(\w^{(k)}|\w^{(k)})&=\sum_{i=1}^{p}\zeta(|w^{(k)}_i|)+\zeta^{\prime}(|w^{(k)}_i|)(|w^{(k+1)}_i|-|w^{(k)}_i|)\\
&-\sum_{i=1}^{p}\zeta(|w^{(k)}_i|)+\zeta^{\prime}(|w^{(k)}_i|)(|w^{(k)}_i|-|w^{(k)}_i|)\\
&=\sum_{i=1}^{p}\langle\zeta^{\prime}(|w^{(k)}_i|)(|w^{(k+1)}_i|-|w^{(k)}_i|)\\
&\geq \sum_{i=1}^{p}\zeta(|w^{(k+1)}_i|)-\zeta(|w^{(k)}_i|)\\
&=r(\w^{(k+1)})-r(\w^{(k)})
\end{split}
\]
Combining the above three inequalities, we have\[
f(\w^{(k)})+r(\w^{(k)})-f(\w^{(k+1)})-r(\w^{(k+1)})\geq \frac{\gamma}{2}\|\w^{(k+1)}-\w^{(k)}\|^2,
\]
which implies that
\[
F(\w^{(k)})-F(\w^{(k+1)})\geq\frac{\gamma}{2}\|\w^{(k+1)}-\w^{(k)}\|^2.
\]
\end{proof}

\subsection{The proof of Lemma \ref{lem:bound-sub}}
\begin{proof}
\begin{enumerate}
\renewcommand{\labelenumi}{(\theenumi)}
\item[(i)]
Recall that
\[
\w^{(k+1)}=\argmin_{\w} \bigg\{Q_f(\w|\w^{(k)})+r(\w)\bigg\}.
\]
Writing down the optimality condition, we have
\begin{equation}
\0 = \nabla Q_f(\w^{(k+1)}|\w^{(k)})+\u^{(k+1)},
\end{equation}
where $\u^{(k+1)}\in \partial r(\w^{(k+1)})$.
Let's rewrite it as follow\[
\nabla Q_f(\w^{(k+1)}|\w^{(k)})-\nabla f(\w^{(k+1)}) + \nabla f(\w^{(k+1)})  + \u^{(k+1)}=\0.
\]
Because $A^{(k+1)}=\nabla f(\w^{(k+1)})-\nabla Q_f(\w^{(k+1)}|\w^{(k)})$, we immediately have
\[
A^{(k+1)}=\nabla f(\w^{(k+1)})+\u^{(k+1)} \in \partial F(\w^{(k+1)}).
\]
\item[(ii)]
With the Lipschitz continuous of $\nabla Q_f(\w|\w^{(k)})$ and $\nabla f(\w)$, we have
\begin{equation}
\left\{
\begin{array}{l}
\|\nabla Q_f(\w^{(k+1)}|\w^{(k)})-\nabla Q_f(\w^{(k)}|\w^{(k)})\|\leq L_{Q_f}\|\w^{(k+1)}-\w^{(k)}\|,\\

\|\nabla f(\w^{(k+1)})-\nabla f(\w^{(k)})\|\leq L_f\|\w^{(k+1)}-\w^{(k)}\|.
\end{array}
\right.
\end{equation}
Hence,
\begin{equation}
\begin{split}
\|\A^{(k+1)}\|&=\|\nabla f(\w^{(k+1)})-\nabla Q_f(\w^{(k+1)}|\w^{(k)})\|\\
&= \|\nabla f(\w^{(k+1)}) -\nabla f(\w^{(k)})+\nabla Q_f(\w^{(k)}|\w^{(k)})-\nabla Q_f(\w^{(k+1)}|\w^{(k)})\|\\
&\leq \|\nabla f(\w^{(k+1)}) -\nabla f(\w^{(k)})\|+\|\nabla Q_f(\w^{(k)}|\w^{(k)})-\nabla Q_f(\w^{(k+1)}|\w^{(k)})\|\\
&\leq L_f\|\w^{(k+1)}-\w^{(k)}\|+L_{Q_f}\|\w^{(k+1)}-\w^{(k)}\|\\
&=(L_{Q_f}+L_f)\|\w^{(k+1)}-\w^{(k)}\|
\end{split}
\end{equation}
\end{enumerate}
\end{proof}

\subsection{The proof of Lemma \ref{lem:share}}
\begin{proof}
\begin{enumerate}
\item[(i)]
By the concavity of $\zeta$, we have
\[
\zeta(|w_i|)\leq \zeta(|w^{(k)}_i|)+\zeta^{\prime}(|w_i^{(k)}|)( |w_i|-|w^{(k)}_i|),
\]
for any $i\in[p]$. We immediately obtain
\[
Q_r(\w|\w^{(k)})\geq r(\w) ~\mbox{and} ~ Q_r(\w^{(k)}|\w^{(k)})= r(\w^{(k)}).
\]
\item[(ii)]
Notice the fact that the subdifferential calculus for separable functions yields the follows \citep{rockafellar1998variational}.
\[
\partial \Big(\sum^p_{i=1} \zeta(w_i)\Big)=\partial \zeta(w_1) \times \partial \zeta(w_2)\cdots \times \partial \zeta(w_n).
\]
Since $r(\w)$ and $Q_r(\w)$ are separable, we consider each dimension independently.
For any $i\in[p]$, if $w_i>0$, we have
\[
\partial_i r(\w^{(k)})=\{\zeta^{\prime}(w_i)\}, ~\partial_i Q_r(\w^{(k)})=\{\zeta^{\prime}(w_i)\}.
\]
Similarly, if $w_i<0$, we have
\[
\partial_i r(\w^{(k)})=\{-\zeta^{\prime}(-w_i)\}, ~\partial_i Q_r(\w^{(k)})=\{-\zeta^{\prime}(-w_i)\}.
\]
For the case $w_i=0$, we have
\[
\partial_i r(\w^{(k)})=[-\zeta^{\prime}(0),\zeta^{\prime}(0)], ~\partial_i Q_r(\w^{(k)})=[-\zeta^{\prime}(0),\zeta^{\prime}(0)].
\]
So
\[
\partial Q_r(\w^{(k)}|\w^{(k)})=\partial r(\w^{(k)}).
\]
\end{enumerate}
\end{proof}

\subsection{The proof of Lemma \ref{cor:bound-sub}}
\begin{proof}
\begin{enumerate}
\item[(i)]
By using
\[
\w^{(k+1)}=\argmin_{\w} \bigg\{Q_f(\w|\w^{(k)})+Q_r(\w|\w^{(k)})\bigg\}.
\]
Alternatively
\[
\w^{(k+1)}=\argmin_{\w}\bigg\{ Q_f(\w|\w^{(k)})+\sum_{i=1}^p\zeta(|w^{(k)}_i|)+\zeta^{\prime}(|w^{(k)}_i|)(|w_i|-|w^{(k)}_i|)\bigg\}.
\]
For the notation, we let $\mu_i^{(k)}=\zeta^{\prime}(|w^{(k)}_i|)\sgn(w^{(k+1)}_i)$.
Then the optimal condition yields
\begin{equation}\label{eq:Opt_condition}
\0=\nabla Q_f(\w^{(k+1)}|\w^{(k)})+(\mu_1^{(k)},\mu_2^{(k)},\cdots,\mu_p^{(k)})^T.
\end{equation}
We rewrite it as
\[
\begin{split}
\0&=\nabla Q_f(\w^{(k+1)}|\w^{(k)})-\nabla f(\w^{(k+1)})+\nabla f(\w^{(k+1)})\\&+(\zeta^{\prime}(|w^{(k)}_1|)\sgn(w^{(k+1)}_1),\zeta^{\prime}(|w^{(k)}_2|)\sgn(w^{(k+1)}_2),\cdots,\zeta^{\prime}(|w^{(k)}_p|)\sgn(w^{(k+1)}_p))^T
\end{split}
\]
On the other hand
\[
b^{(k)}_i=\sgn(w_i^{(k+1)})(\zeta^{\prime}(|w_i^{(k)}|)-\zeta^{\prime}(|w_i^{(k+1)}|),
\]
for $i\in[p]$.
So
\[
\begin{split}
\0&=\nabla f(\w^{(k+1)})+(\zeta^{\prime}(|w^{(k+1)}_1|)\sgn(w^{(k+1)}_1),\zeta^{\prime}(|w^{(k+1)}_2|)\sgn(w^{(k+1)}_2),\cdots,\zeta^{\prime}(|w^{(k+1)}_p|)\sgn(w^{(k+1)}_p))^T\\&+\nabla Q_f(\w^{(k+1)}|\w^{(k)})-\nabla f(\w^{(k+1)})+(b^{(k)}_1,b^{(k)}_2,\cdots,b^{(k)}_p)^T.
\end{split}
\]
Then, we have
\[
\begin{split}
\nabla f(\w^{(k+1)})&+(\zeta^{\prime}(|w^{(k+1)}_1|)\sgn(w^{(k+1)}_1),\zeta^{\prime}(|w^{(k+1)}_2|)\sgn(w^{(k+1)}_2),\cdots,\zeta^{\prime}(|w^{(k+1)}_p|)\sgn(w^{(k+1)}_p))^T\\
&=\nabla f(\w^{(k+1)})-\nabla Q_f(\w^{(k+1)}| \w^{(k)})- (b^{(k)}_1,b^{(k)}_2,\cdots,b^{(k)}_p)^T\\
&=B^{(k+1)}
\end{split}
\]
Notice that
\[
\begin{split}
\nabla f(\w^{(k+1)})&+(\zeta^{\prime}(|w^{(k+1)}_1|)\sgn(w^{(k+1)}_1),\zeta^{\prime}(|w^{(k+1)}_2|)\sgn(w^{(k+1)}_2),\cdots,\zeta^{\prime}(|w^{(k+1)}_p|)\sgn(w^{(k+1)}_p))^T\\&\in\partial F(\w^{(k+1)})
\end{split}
\]
So we have
\[
B^{(k+1)}\in \partial F(\w^{(k+1)}).
\]
\item[(ii)]
Similarly, with the Lipschitz continuous of $\nabla Q_f(\w|\w^{(k)})$ , $\nabla f(\w)$ and $\zeta^{\prime}(t)(t\geq0)$ we have 

\begin{equation}
\left\{
\begin{array}{l}
\|\nabla Q_f(\w^{(k+1)}|\w^{(k)})-\nabla Q_f(\w^{(k)}|\w^{(k)})\|\leq L_{Q_f}\|\w^{(k+1)}-\w^{(k)}\|,\\

\|\nabla f(\w^{(k+1)})-\nabla f(\w^{(k)})\|\leq L_f\|\w^{(k+1)}-\w^{(k)}\|,\\

|\zeta^{\prime}(t_1)-\zeta^{\prime}(t_2)|\leq L_{\zeta}|t_1-t_2|.
\end{array}
\right.
\end{equation}
Now, we are ready to bound the subgradient $B^{(k+1)}$
\begin{equation}\label{eq:BoudB}
\begin{split}
\|B^{(k+1)}\|&=\|\nabla f(\w^{(k+1)})-\nabla Q_f(\w^{(k+1)}| \w^{(k)})- (b^{(k)}_1,b^{(k)}_2,\cdots,b^{(k)}_p)^T\|\\
&=\|\nabla f(\w^{(k+1)})-\nabla f(\w^{(k)})+\nabla Q_f(\w^{(k)}|\w^{(k)})\\&-\nabla Q_f(\w^{(k+1)}| \w^{(k)})- (b^{(k)}_1,b^{(k)}_2,\cdots,b^{(k)}_p)^T\|\\
&\leq\|\nabla f(\w^{(k+1)})-\nabla f(\w^{(k)})\|+\|\nabla Q_f(\w^{(k)}|\w^{(k)})-\nabla Q_f(\w^{(k+1)}|\w^{(k)})\|\\
&+\|(b^{(k)}_1,b^{(k)}_2,\cdots,b^{(k)}_p)^T\|\\
&\leq L_f\|\w^{(k+1)}-\w^{(k)}\|+L_{Q_f}\|\w^{(k+1)}-\w^{(k)}\|+\|(b^{(k)}_1,b^{(k)}_2,\cdots,b^{(k)}_p)^T\|\\
\end{split}
\end{equation}
Then let's bound $\|(b^{(k)}_1,b^{(k)}_2,\cdots,b^{(k)}_p)^T\|$.
For each $i\in[p]$, we have
\[
b^{(k)}_i=\sgn(w_i^{(k+1)})(\zeta^{\prime}(|w_i^{(k)}|)-\zeta^{\prime}(|w_i^{(k+1)}|)).
\]
By using $|\sgn(w_i^{(k+1)})|\leq1$, we have
\[
\begin{split}
|b^{(k)}_i|&=|\sgn(w_i^{(k+1)})(\zeta^{\prime}(|w_i^{(k)}|)-\zeta^{\prime}(|w_i^{(k+1)}|))|\\
&\leq|(\zeta^{\prime}(|w_i^{(k)}|)-\zeta^{\prime}(|w_i^{(k+1)}|))|\\
&\leq L_\zeta||w^{(k+1)}_i|-|w^{(k)}_i||
\end{split}
\]
Then
\begin{equation}\label{eq:Boundb}
\begin{split}
\|(b^{(k)}_1,b^{(k)}_2,\cdots,b^{(k)}_p)^T\|&\leq L_\zeta\|(|w_1^{(k+1)}|-|w_1^{(k)}|,|w_2^{(k+1)}|-|w_2^{(k)}|,\cdots,|w_p^{(k+1)}|-|w_p^{(k)}|)^T\|\\
&\leq L_\zeta\|(|w_1^{(k+1)}-w_1^{(k)}|,|w_2^{(k+1)}-w_2^{(k)}|,\cdots,|w_p^{(k+1)}-w_p^{(k)}|)^T\| \\
&=L_\zeta\|(w_1^{(k+1)}-w_1^{(k)},w_2^{(k+1)}-w_2^{(k)},\cdots,w_p^{(k+1)}-w_p^{(k)})^T\|\\
&=L_\zeta\|\w^{(k+1)}-\w^{(k)}\|
\end{split}
\end{equation}
Combining (\ref{eq:BoudB}) and (\ref{eq:Boundb}), 
\[
\|B^{(k+1)}\|\leq (L_{Q_f}+L_f+L_\zeta)\|\w^{(k+1)}-\w^{(k)}\|.
\]
\end{enumerate}
\end{proof}

\subsection{The proof of Lemma \ref{lemm:crictical}}
\begin{proof}
Since the F(\w) is coercive, the sequence $\left\{\w^{(k)}\right\}_{k\in\mathbb N}$ is bounded. Therefore there exists an increasing sequence $\left\{n_k\right\}_{k\in\mathbb N}$ such that
\[
\lim_{k\to \infty}\w^{(n_k)}=\w^*.
\]
Recall that $F(\w)=f(\w)+\lambda\sum^p_{i=1} \zeta(w_{i})$ is continuous.
We have
\[
\lim_{k\to\infty} F(\w^{(n_k)})= F(\w^*).
\]
On the other hand, we know $A^{(k)}\in\partial F(\w^{(k)}), B^{(k)}\in\partial F(\w^{(k)})$. Moreover, from Lemma (\ref{lem:bound-sub}) and Lemma (\ref{cor:bound-sub}) it can be seen that as $k \to \infty$, $A^{(k)}\to \0$ and $B^{(k)} \to \0$.
Remember that $\partial F$ is close. So $\0 \in \partial F(\w^*)$, which contributes to that $\w^{*}$ is a critical point of F.
\end{proof}

\subsection{Uniformized K$\L{}$ property}
Before providing the global convergence result, 
we first introduce a class of concave and continuous functions. Let $\eta\in(0,+\infty]$. We are concerned with $\Phi_\eta$ which contain the class of all concave and continuous functions $\phi:[0,\eta)\to \RB_+$  satisfying the following properties:

\begin{enumerate}
\item[(a)] $\phi(0)=0$ and continuous at 0; 
\item[(b)] $\phi$ is $C^1$ on $(0,\eta)$; 
\item[(c)] $\phi^{\prime}(t)>0$, $\forall$ t $\in$ (0,$\eta$).
\end{enumerate}

\begin{lemma}[\cite{bolte2013proximal}]\label{uniKL}
Suppose $\Omega$ is a compact set and let $F:\RB^p\to(-\infty,\infty]$ be a lower semi-continuous function. Moreover, $F$ is constant on $\Omega$ and satisfy KL property at each point of $\Omega$. Then there exist $\epsilon>0,\eta>0$ and $\phi\in\Phi_\eta$ such that for all $\bar \u$ in $\Omega$ and all $\u$ in 
\[
\bigg\{\u \in \RB^p:~dist(\u,\Omega)<\epsilon\bigg\}\bigcap\bigg\{\u:F(\bar \u)<F(\u)<F(\bar \u)+\eta\bigg\}
\]
one has,
\begin{eqnarray}\label{lemm:con}
\phi^{\prime}(F(\u)-F(\bar \u)) dist(\0,\partial F(\u)) \geq 1.
\end{eqnarray}
\end{lemma}

\subsection{The proof of Theorem \ref{GlobalCon}}
\begin{proof}\label{pf:main}
As is known, there exists an increasing sequence $\left\{n_{k}\right\}_{k\in\mathbb N}$ such that $\left\{\w^{(n_k)}\right\}$ converges to $\w^*$.
Suppose that there exists an integer $n_0$ satisfy that $F(\w^{(n_0)})=F(\w^*)$. Then it is clear that for any integer $N>n_0$, $F(\w^{(N)})=F(\w^*)$ holds. Then it is trivial to achieve the convergent sequence. Otherwise, we consider the case that $F(\w^{(k)})>F(\w^*)$, $\forall k \in \mathbb N$. Because the sequence $\{F(\w^{(k)})\}_{k\in\mathbb N}$ is convergent, it is clear that for any $\eta>0$, there exist one integer $m$ such that $F(\w^{(k)})<F(\w^*)+\eta$ for all $k>m$. By using lemma \ref{lemm:crictical}, we have $\lim_{k\to\infty}dist(\w^{(k)},\mathcal M(\w^{(0)}))=0$ which implies that for any $\epsilon>0$ there exists a positive integer $n$ such that $dist(\w^{(k)},\mathcal M(\w^{(0)}))<\epsilon$ for all $k>n$. Let $l=\max\left\{m,n\right\}$.
Then for any $k>l$, we have
\[
\phi^{\prime}\Big(F(\w^{(k)})-F(\w^*)\Big)dist(\0,\partial F(\w^{(k)}))\geq 1.
\]
By the Lemma \ref{lem:bound-sub} and Lemma \ref{cor:bound-sub}, we have
\begin{equation}\label{eq:lb}
\phi^{\prime}\Big(F(\w^{(k)})-F(\w^*)\Big) \geq \rho \Big\|\w^{(k)}-\w^{(k-1)}\Big\|^{-1},
\end{equation}
where $\rho=\frac{1}{L_{Q_f}+L_f+L_\zeta}$.

we let~ $d_{k,k+1}$$\denote\phi(F(\w^{(k)})-F(\w^*))-\phi(F(\w^{(k+1)})-F(\w^*))$. With the property of concave functions, we have
\begin{equation}\label{eq:concave}
\begin{split}
d_{k,k+1}&\geq\phi^{\prime}(F(\w^{(k)})-F(\w^*))(F(\w^{k})-F(\w^{(k+1)})) \\
         &\geq \frac{\gamma}{2}\phi^{\prime}(F(\w^{(k)})-F(\w^*))\|\w^{(k+1)}-\w^{(k)}\|^2 \\
         &\geq \frac{\gamma}{2(L_f+L_{Q_f}+L_\zeta)}\|\w^{(k)}-\w^{(k-1)}\|^{-1}\|\w^{(k+1)}-\w^{(k)}\|^2
\end{split}
\end{equation}
That is
\[
  Md_{k,k+1}\|\w^{(k)}-\w^{k-1}\|\geq\|\w^{(k+1)}-\w^{(k)}\|^2
\]
where $M$ = $\frac{2(L_f+L_{Q_f}+L_\zeta)}{\gamma}$.
Notice that
\[
  Md_{k,k+1}\|\w^{(k)}-\w^{(k-1)}\|\leq\Big(\frac{Md_{k,k+1}+\|\w^{(k)}-\w^{(k-1)}\|}{2}\Big)^2.
\]
So we have
\[
2\|\w^{(k+1)}-\w^{(k)}\|\leq Md_{k,k+1}+\|\w^{(k)}-\w^{(k-1)}\|.
\]
Then
\[
\begin{split}
\sum^{\infty}_{k=l+1}\|\w^{(k+1)}-\w^{(k)}\|&\leq M\sum^{\infty}_{k=l+1}d_{k,k+1}+\sum^{\infty}_{k=l+1}\Big(\|\w^{(k)}-\w^{(k-1)}\|-\|\w^{(k+1)}-\w^{(k)}\|\Big) \\
&\leq Md_{l+1,\infty}+\|\w^{(l+1)}-\w^{(l)}\| \\
&\leq M\phi\Big(F(\w^{(l+1)})-F(\w^*)\Big)+\|\w^{(l+1)}-\w^{(l)}\|
\end{split}
\]
Let $l \to \infty$. Since $\lim_{l \to \infty}\|\w^{(l+1)}-\w^{(l)}\|=0$ and $\lim_{l \to \infty}F(\w^{(l+1)})=F(\w^*)$, it is clear that
\[
\lim_{l \to \infty}\sum^{\infty}_{k=l+1}\|\w^{(k+1)}-\w^{(k)}\|=0
\]
So
\[
 \sum^{\infty}_{k=0}\|\w^{(k+1)}-\w^{(k)}\|<\infty.
\]
Then, we have 
\[
\lim_{m\to\infty} \sum^l_{k=m}\|\w^{(k+1)}-\w^{(k)}\|=0,
\]
for any $m<l$.
This suggests that $\left\{\w^{(k)}\right\}_{k\in \mathbb N}$ is Cauchy sequence. As a result, it is a convergent sequence that converges to $\w^*$.
\end{proof}

\subsection{The proof of Theorem \ref{theorem:convergenceRate}}
This is a classical result of K\L{} function.
Since the corresponding function 
\[
\phi(t) = c t^{1-\theta}, \theta \in[0, 1).
\]

As \cite{attouch2009convergence}, the conclusions of Theorem \ref{theorem:convergenceRate} hold.

\subsection{The proof of Theorem \ref{the:cccp}}
\begin{proof}
Recall that
\[
\w^{(k+1)} = \; \argmin_{\w}~~~~\bigg\{ \delta_{\mathcal C}(\w)+u(\w)-\nabla v(\w^{(k)})^T\w \bigg\}.
\]
We immediately have\[
\0=\g^{(k+1)}+\nabla u(\w^{(k+1)})-\nabla v(\w^{(k)}),
\]
where $\g^{(k+1)}\in \partial \delta_{\mathcal C}(\w^{(k+1)})$.

Because $\delta_{\mathcal C}(\w)+u(\w)-\nabla v(\w^{(k)})^T\w$ is $\gamma$-strongly convex, we have
\begin{equation}
\begin{split}
\delta_{\mathcal C}(\w^{(k)})+u(&\w^{(k)})-\nabla v(\w^{(k)})^T\w^{(k)}-\delta_{\mathcal C}(\w^{(k+1)})-u(\w^{(k+1)})+\nabla v(\w^{(k)})^T\w^{(k+1)}\\
&\geq\langle\0,\w^{(k)}-\w^{(k+1)}\rangle+\frac{\gamma}{2}\|\w^{(k+1)}-\w^{(k)}\|^2.
\end{split}
\end{equation}
By the convexity of $v(\w)$, we obtain
\[
\nabla v(\w^{(k)})^T\w^{(k+1)}-\nabla v(\w^{(k)})^T\w^{(k)}\leq v(\w^{(k+1)})-v(\w^{(k)}).
\]
Thus,
\begin{equation}\label{eq:cccp_descent}
\begin{split}
\Big(\underbrace{\delta_{\mathcal C}(\w^{(k)})+u(\w^{(k)})-v(\w^{(k)})}_{F(\w^{(k)})}\Big)&-\Big(\underbrace{\delta_{\mathcal C}(\w^{(k+1)})+u(\w^{(k+1)})-v(\w^{(k+1)}}_{F(\w^{(k+1)})}\Big)\geq \frac{\gamma}{2}\|\w^{(k+1)}-\w^{(k)}\|^2.
\end{split}
\end{equation}
On the other hand, we known that 
\[
\0=\g^{(k+1)}+\nabla u(\w^{(k+1)})-\nabla v(\w^{(k+1)}) + \nabla v(\w^{(k+1)})-\nabla v(\w^{(k)}).
\]
Let's denote $\nabla v(\w^{(k)})-\nabla v(\w^{(k+1)}$ as $C^{(k+1)}$. Then
$C^{(k+1)}\in\partial F(\w^{(k+1)}).$
\begin{equation}\label{eq:cccp_subbound}
\begin{split}
\|C^{(k+1)}\|&=\|\nabla v(\w^{(k)})-\nabla v(\w^{(k+1)})\| \leq L_v\|\w^{(k+1)}-\w^{(k)}\|
\end{split}
\end{equation}
Next we prove $\mathcal M(\w^{(0)})$ are the subset of the crit $F(\w)$.
Because of the coerciveness of the Function $F(\w)$, there exists a bounded sequence $\{\w^{n_k}\}_{k\in\mathbb N}$, which satisfies that $\lim_{k\to\infty}\w^{n_{k}}=\bar\w$.
Since $\delta_{\mathcal C}(\w)$ is lower semicontinuous, we have
\begin{equation}\label{eq:liminf}
\liminf_{k\to\infty} \delta_{\mathcal C}(\w^{n_k})\geq \delta_{\mathcal C}(\bar\w)
\end{equation}
On the other hand,
\[
\delta_{\mathcal C}(\w^{(k+1)})+u(\w^{(k+1)})-\nabla v(\w^{(k)})^T\w^{(k+1)}\leq
\delta_{\mathcal C}(\bar\w)+u(\bar\w)-\nabla v(\w^{(k)})^T\bar\w
\]
Rewrite the above formulation, we obtain
\[
\delta_{\mathcal C}(\w^{(k+1)})\leq
\delta_{\mathcal C}(\bar\w)-(u(\w^{(k+1)})-u(\bar\w))+\nabla v(\w^{(k)})^T(\w^{(k+1)}-\bar\w).
\]
Substitute $k$ with $n_{k}-1$. By the fact $u,v$ are $C^1$ functions, we have
\begin{equation}\label{eq:limsup}
\limsup_{k\to\infty}\delta(\w^{n_k})\leq \delta_{\mathcal C}(\bar\w).
\end{equation}
Combing (\ref{eq:liminf}) and (\ref{eq:limsup}), we immediately have
\[
\lim_{k\to\infty} \delta_{\mathcal C}(\w^{n_k})=\delta_{\mathcal C}(\bar{\w}).
\]
Notice that $u(\w),v(\w)$ are continuous, we have
\[
\lim_{k\to\infty} F(\w^{n_k})=F(\bar\w).
\]
(\ref{eq:cccp_subbound}) implies that $C^{(k)}\to\0$ as $k\to\infty$. Moreover, $C^{(k)}\in\partial F(\w^{(k)})$. Remember the closeness of $\partial F(\w)$, we have $\0\in\partial F(\bar{\w})$. So $\mathcal M(\w^{(0)})$ are the subset of the crit $F(\w)$.
With (\ref{eq:cccp_descent}) and (\ref{eq:cccp_subbound}) ready, the next proof is the same as that of  Theorem \ref{GlobalCon}.
\end{proof}


\vskip 0.2in
\bibliography{GNMM}

\begin{thebibliography}{37}
\providecommand{\natexlab}[1]{#1}
\providecommand{\url}[1]{\texttt{#1}}
\expandafter\ifx\csname urlstyle\endcsname\relax
  \providecommand{\doi}[1]{doi: #1}\else
  \providecommand{\doi}{doi: \begingroup \urlstyle{rm}\Url}\fi

\bibitem[Armagan et~al.(2013)Armagan, Dunson, and Lee]{armagan2013generalized}
Artin Armagan, David~B Dunson, and Jaeyong Lee.
\newblock Generalized double pareto shrinkage.
\newblock \emph{Statistica Sinica}, 23\penalty0 (1):\penalty0 119, 2013.

\bibitem[Attouch and Bolte(2009)]{attouch2009convergence}
Hedy Attouch and J{\'e}r{\^o}me Bolte.
\newblock On the convergence of the proximal algorithm for nonsmooth functions
  involving analytic features.
\newblock \emph{Mathematical Programming}, 116\penalty0 (1-2):\penalty0 5--16,
  2009.

\bibitem[Attouch et~al.(2010)Attouch, Bolte, Redont, and
  Soubeyran]{attouch2010proximal}
H{\'e}dy Attouch, J{\'e}r{\^o}me Bolte, Patrick Redont, and Antoine Soubeyran.
\newblock Proximal alternating minimization and projection methods for
  nonconvex problems: an approach based on the kurdyka-lojasiewicz inequality.
\newblock \emph{Mathematics of Operations Research}, 35\penalty0 (2):\penalty0
  438--457, 2010.

\bibitem[Beck and Teboulle(2009)]{beck2009fast}
Amir Beck and Marc Teboulle.
\newblock A fast iterative shrinkage-thresholding algorithm for linear inverse
  problems.
\newblock \emph{SIAM Journal on Imaging Sciences}, 2\penalty0 (1):\penalty0
  183--202, 2009.

\bibitem[Bochnak et~al.(1998)Bochnak, Coste, and Roy]{bochnak1998real}
Jacek Bochnak, Michel Coste, and Marie-Fran{\c{c}}oise Roy.
\newblock \emph{Real algebraic geometry}.
\newblock Springer, 1998.

\bibitem[Bolte et~al.(2007)Bolte, Daniilidis, and Lewis]{bolte2007lojasiewicz}
J{\'e}r{\^o}me Bolte, Aris Daniilidis, and Adrian Lewis.
\newblock The lojasiewicz inequality for nonsmooth subanalytic functions with
  applications to subgradient dynamical systems.
\newblock \emph{SIAM Journal on Optimization}, 17\penalty0 (4):\penalty0
  1205--1223, 2007.

\bibitem[Bolte et~al.(2013)Bolte, Sabach, and Teboulle]{bolte2013proximal}
J{\'e}r{\^o}me Bolte, Shoham Sabach, and Marc Teboulle.
\newblock Proximal alternating linearized minimization for nonconvex and
  nonsmooth problems.
\newblock \emph{Mathematical Programming}, pages 1--36, 2013.

\bibitem[Candes et~al.(2008)Candes, Wakin, and Boyd]{candes2008enhancing}
Emmanuel~J Candes, Michael~B Wakin, and Stephen~P Boyd.
\newblock Enhancing sparsity by reweighted ℓ 1 minimization.
\newblock \emph{Journal of Fourier analysis and applications}, 14\penalty0
  (5-6):\penalty0 877--905, 2008.

\bibitem[Chartrand and Yin(2008)]{chartrand2008iteratively}
Rick Chartrand and Wotao Yin.
\newblock Iteratively reweighted algorithms for compressive sensing.
\newblock In \emph{Acoustics, speech and signal processing, 2008. ICASSP 2008.
  IEEE international conference on}, pages 3869--3872. IEEE, 2008.

\bibitem[Combettes and Pesquet(2011)]{combettes2011proximal}
Patrick~L Combettes and Jean-Christophe Pesquet.
\newblock Proximal splitting methods in signal processing.
\newblock In \emph{Fixed-point algorithms for inverse problems in science and
  engineering}, pages 185--212. Springer, 2011.

\bibitem[Fan and Li(2001)]{fan2001variable}
Jianqing Fan and Runze Li.
\newblock Variable selection via nonconcave penalized likelihood and its oracle
  properties.
\newblock \emph{Journal of the American Statistical Association}, 96\penalty0
  (456):\penalty0 1348--1360, 2001.

\bibitem[Fukushima and Mine(1981)]{fukushima1981generalized}
Masao Fukushima and Hisashi Mine.
\newblock A generalized proximal point algorithm for certain non-convex
  minimization problems.
\newblock \emph{International Journal of Systems Science}, 12\penalty0
  (8):\penalty0 989--1000, 1981.

\bibitem[Gasso et~al.(2009)Gasso, Rakotomamonjy, and Canu]{gasso2009recovering}
Gilles Gasso, Alain Rakotomamonjy, and St{\'e}phane Canu.
\newblock Recovering sparse signals with a certain family of nonconvex
  penalties and dc programming.
\newblock \emph{Signal Processing, IEEE Transactions on}, 57\penalty0
  (12):\penalty0 4686--4698, 2009.

\bibitem[Gong et~al.(2013)Gong, Zhang, Lu, Huang, and Ye]{gong2013general}
Pinghua Gong, Changshui Zhang, Zhaosong Lu, Jianhua Huang, and Jieping Ye.
\newblock A general iterative shrinkage and thresholding algorithm for
  non-convex regularized optimization problems.
\newblock In \emph{Proceedings of The 30th International Conference on Machine
  Learning}, pages 37--45, 2013.

\bibitem[Iusem(1999)]{iusem1999augmented}
AN~Iusem.
\newblock Augmented lagrangian methods and proximal point methods for convex
  optimization.
\newblock \emph{Investigaci{\'o}n Operativa}, 8:\penalty0 11--49, 1999.

\bibitem[Kurdyka(1998)]{kurdyka1998gradients}
Krzysztof Kurdyka.
\newblock On gradients of functions definable in o-minimal structures.
\newblock In \emph{Annales de l'institut Fourier}, volume~48, pages 769--783.
  Institut Fourier, 1998.

\bibitem[Lanckriet and Sriperumbudur(2009)]{lanckriet2009convergence}
Gert~R Lanckriet and Bharath~K Sriperumbudur.
\newblock On the convergence of the concave-convex procedure.
\newblock In \emph{Advances in neural information processing systems}, pages
  1759--1767, 2009.

\bibitem[Lange(2004)]{lange2004optimization}
Kenneth Lange.
\newblock Optimization. springer texts in statistics.
\newblock 2004.

\bibitem[Lange et~al.(2000)Lange, Hunter, and Yang]{lange2000optimization}
Kenneth Lange, David~R Hunter, and Ilsoon Yang.
\newblock Optimization transfer using surrogate objective functions.
\newblock \emph{Journal of computational and graphical statistics}, 9\penalty0
  (1):\penalty0 1--20, 2000.

\bibitem[Lemaire(1989)]{lemaire1989proximal}
Bernard Lemaire.
\newblock The proximal algorithm.
\newblock \emph{International series of numerical mathematics}, 87:\penalty0
  73--87, 1989.

\bibitem[Lewis and Wright(2008)]{lewis2008proximal}
Adrian~S Lewis and Stephen~J Wright.
\newblock A proximal method for composite minimization.
\newblock \emph{arXiv preprint arXiv:0812.0423}, 2008.

\bibitem[{\L}ojasiewicz(1993)]{lojasiewicz1993geometrie}
Stanislas {\L}ojasiewicz.
\newblock Sur la g{\'e}om{\'e}trie semi-et sous-analytique.
\newblock In \emph{Annales de l'institut Fourier}, volume~43, pages 1575--1595.
  Institut Fourier, 1993.

\bibitem[Lojasiewicz(1963)]{lojasiewicz1963propriete}
Stanislaw Lojasiewicz.
\newblock Une propri{\'e}t{\'e} topologique des sous-ensembles analytiques
  r{\'e}els.
\newblock \emph{Colloques du CNRS, Les {\'e}quations aux d{\'e}riv{\'e}s
  partielles}, 117, 1963.

\bibitem[Mairal(2013)]{mairal2013optimization}
Julien Mairal.
\newblock Optimization with first-order surrogate functions.
\newblock In \emph{ICML 2013-International Conference on Machine Learning},
  volume~28, pages 783--791, 2013.

\bibitem[Mazumder et~al.(2011)Mazumder, Friedman, and
  Hastie]{mazumder2011sparsenet}
Rahul Mazumder, Jerome~H Friedman, and Trevor Hastie.
\newblock Sparsenet: Coordinate descent with nonconvex penalties.
\newblock \emph{Journal of the American Statistical Association}, 106\penalty0
  (495), 2011.

\bibitem[Nesterov and Nesterov(2004)]{nesterov2004introductory}
Yurii Nesterov and I͡U~E Nesterov.
\newblock \emph{Introductory lectures on convex optimization: A basic course},
  volume~87.
\newblock Springer, 2004.

\bibitem[Parikh and Boyd(2013)]{parikh2013proximal}
Neal Parikh and Stephen Boyd.
\newblock Proximal algorithms.
\newblock \emph{Foundations and Trends in Optimization}, 1\penalty0
  (3):\penalty0 123--231, 2013.

\bibitem[Rockafellar(1976)]{rockafellar1976monotone}
R~Tyrrell Rockafellar.
\newblock Monotone operators and the proximal point algorithm.
\newblock \emph{SIAM journal on control and optimization}, 14\penalty0
  (5):\penalty0 877--898, 1976.

\bibitem[Rockafellar et~al.(1998)Rockafellar, Wets, and
  Wets]{rockafellar1998variational}
R~Tyrrell Rockafellar, Roger J-B Wets, and Maria Wets.
\newblock \emph{Variational analysis}, volume 317.
\newblock Springer, 1998.

\bibitem[Vaida(2005)]{vaida2005parameter}
Florin Vaida.
\newblock Parameter convergence for em and mm algorithms.
\newblock \emph{Statistica Sinica}, 15\penalty0 (3):\penalty0 831, 2005.

\bibitem[Wu(1983)]{wu1983convergence}
CF~Jeff Wu.
\newblock On the convergence properties of the em algorithm.
\newblock \emph{The Annals of statistics}, pages 95--103, 1983.

\bibitem[Xu and Yin(2013)]{xu2013block}
Yangyang Xu and Wotao Yin.
\newblock A block coordinate descent method for regularized multiconvex
  optimization with applications to nonnegative tensor factorization and
  completion.
\newblock \emph{SIAM Journal on imaging sciences}, 6\penalty0 (3):\penalty0
  1758--1789, 2013.

\bibitem[Yuille and Rangarajan(2003)]{yuille2003concave}
Alan~L Yuille and Anand Rangarajan.
\newblock The concave-convex procedure.
\newblock \emph{Neural Computation}, 15\penalty0 (4):\penalty0 915--936, 2003.

\bibitem[Zhang(2010{\natexlab{a}})]{zhang2010nearly}
Cun-Hui Zhang.
\newblock Nearly unbiased variable selection under minimax concave penalty.
\newblock \emph{The Annals of Statistics}, pages 894--942, 2010{\natexlab{a}}.

\bibitem[Zhang et~al.(2012)Zhang, Zhang, et~al.]{zhang2012general}
Cun-Hui Zhang, Tong Zhang, et~al.
\newblock A general theory of concave regularization for high-dimensional
  sparse estimation problems.
\newblock \emph{Statistical Science}, 27\penalty0 (4):\penalty0 576--593, 2012.

\bibitem[Zhang(2010{\natexlab{b}})]{zhang2010analysis}
Tong Zhang.
\newblock Analysis of multi-stage convex relaxation for sparse regularization.
\newblock \emph{The Journal of Machine Learning Research}, 11:\penalty0
  1081--1107, 2010{\natexlab{b}}.

\bibitem[Zou and Li(2008)]{zou2008one}
Hui Zou and Runze Li.
\newblock One-step sparse estimates in nonconcave penalized likelihood models.
\newblock \emph{Annals of statistics}, 36\penalty0 (4):\penalty0 1509, 2008.

\end{thebibliography}

\end{document}